\documentclass[11pt]{article}
\usepackage{latexsym,amsfonts,amssymb,amsmath,euscript,amsthm}
\usepackage[applemac]{inputenc}  %para mac
\usepackage[english,french]{babel}  
\usepackage{graphicx} 
\usepackage{subfigure}% inserir imagens
\usepackage{hyperref}         % para acrescentar links
\usepackage{color,fancybox}
\usepackage{bm}
\usepackage{float}

\newcommand{\R} {\mathbb{R}}
\newcommand{\N} {\mathbb{N}}

\newcommand{\D}{\displaystyle}

\newcommand{\pa}{\noindent}

\newcommand{\uK}{K^{\sup}}
\newcommand{\oK}{K^{\inf}}

\newcommand{\xdownarrow}[1]{%
  {\left\downarrow\vbox to #1{}\right.\kern-\nulldelimiterspace}
}

\newcommand{\eproof}{ \vspace{0.3cm} \nolinebreak  \hfill   \raisebox{-\baselineskip}{\llap{\openbox}}}

\numberwithin{equation}{section}
\newtheorem{teo}{Theorem}[section]
\newtheorem{lem}[teo]{Lemma}
\newtheorem{cor}[teo]{Corollary}
\newtheorem{prop}[teo]{Proposition}

\newtheorem{re}[teo]{Remark}
\newtheorem{ex}[teo]{Example}
\newtheorem*{assf*}{Assumption}

\catcode`\@=11
\def\newcodes@{\catcode`\\=12 \catcode`\{=12 \catcode`\}=12 \catcode`\#=12
 \catcode`\%=12\relax}
\def\oldcodes@{\catcode`\\=0 \catcode`\{=1 \catcode`\}=2 \catcode`\#=6
 \catcode`\%=14\relax}
\def\comment{\newcodes@\endlinechar=10 \comment@}
{\lccode`\!=`\\
\lowercase{\gdef\comment@#1^^J{\comment@@#1!endcomment\comment@@@}%
\gdef\comment@@#1!endcomment{\futurelet\next\comment@@@}%
\gdef\comment@@@#1\comment@@@{\ifx\next\comment@@@\let
\next=\comment@\else\def\next{\oldcodes@\endlinechar=`\^^M\relax}%
 \fi\next}}}
\catcode`\@=\active

\title{Boundedness of solutions of nonautonomous  \\ degenerate logistic equations
}

\author{Jos\'{e} M. Arrieta \\Departamento de An\'alisis y Matem\'atica Aplicada, \\ Universidad Complutense de Madrid, 28040 Madrid, Spain \\ and \\ Instituto de Ciencias Matem\'aticas
CSIC-UAM-UC3M-UCM, \\ C/Nicol\'as Cabrera 13-15, Cantoblanco, 28049 Madrid, Spain. \\ e-mail: arrieta at mat.ucm.es  \\  \\ 
Marcos Molina-Rodríguez  \\ Departamento de An\'alisis y Matem\'atica Aplicada \\ Universidad Complutense de Madrid, 28040 Madrid, Spain. \\e-mail: marcos.molina.rodriguez at gmail.com  \\ \\
Lucas A. Santos \\  Instituto Federal da Para\'iba,   \\ Jo\?ao Pessoa PB 58051-900	 PB, Brazil \\ e-mail: lucas.santos at ifpb.edu.br
\\ \\ \\
\rightline {\sl Dedicated to the memory of Genevi\`eve Raugel}}

\date{}

\begin{document}

\maketitle

{\footnotesize 
\par\noindent {\bf Abstract.}
In this work we analyze the boundedness properties of the solutions of a nonautonomous parabolic degenerate logistic equation in a bounded domain. The equation is degenerate in the sense that the  logistic nonlinearity vanishes in a moving region, $K(t)$,  inside the domain.  The boundedness character of the solutions depends not only on, roughly speaking, the first eigenvalue of the Laplace operator in $K(t)$ but also on the way this moving set evolves inside the domain and in particular on the speed at which it moves. 

 \vskip 0.5\baselineskip

\vspace{11pt}

\vspace{6pt}

}

\numberwithin{equation}{section}
\newtheorem{theorem}{Theorem}[section]
\newtheorem{lemma}[theorem]{Lemma}
\newtheorem{corollary}[theorem]{Corollary}
\newtheorem{proposition}[theorem]{Proposition}
\newtheorem{definition}[theorem]{Definition}
\newtheorem{remark}[theorem]{Remark}
\allowdisplaybreaks

\section{Introduction}
\selectlanguage{english}
In this work we study the behavior of the solutions of a nonautonomous evolution problem of the type

\begin{equation}\label{PreEq1}\left\{\begin{array}{lr}
u_t-\Delta u = \lambda u - n(t,x) |u|^{\rho-1}u, & x\in\Omega, \; t>t_0,\\
u=0, & x\in\partial\Omega, \; t>t_0,\\
u(t_0,x)=u_0(x)\geq 0, & x\in \Omega,
\end{array}\right.\end{equation}

\pa where $\Omega$ is a bounded and smooth domain in $\R^N$, for some positive integer $N$, $\rho>1, t_0,\lambda\in\R$, $n\in L^\infty(\R\times\Omega)$ with $n\geq 0$ and $0\leq u_0\in L^q(\Omega)$ where $1<q\leq \infty$ will be chosen appropriately large enough. We consider Dirichlet boundary conditions but the analysis could be developed for other boundary conditions.

The set where the function $n(t,\cdot)$ vanishes is denoted by $K(t)$, that is, 
\begin{equation}\label{Kt}
K(t):=\left\{x\in \overline{\Omega} : n(t,x)=0  \right\}, \; \; t\in \R.
\end{equation}
that we assume is compact for every $t\in \R$, probably empty for some times or even the entire domain for other intervals of time. 

%\pa Therefore, within $K(t)$, equation (\ref{PreEq1}) becomes linear.  
Notice that we are not assuming any special structure of the function $n(t,x)$. We are not restricting ourselves to the periodic case, that is $n(t+T,x)=n(t,x)$ for certain $T>0$ or even to quasiperiodic or almost periodic cases. These are relevant and very interesting cases, but we are aiming to cover general nonatuonomous situation, including these ones.

From a population dynamics interpretation, see \cite{cantrellcosner}, equation \eqref{PreEq1} represents the evolution of the density of certain population, $u(x,t)$, in certain habitat given by the domain $\Omega$. The species diffuses inside the domain and the fact that the boundary condition are homogeneous Dirichlet represents that the exterior of the habitat is very hostile for the species and no individual of the species survives outside the habitat.  The term $\lambda u$ represent the natural growth rate of the species, that we assume to be constant and the term $-n(t,x) |u|^{\rho -1}u$ represents a logistic term that opposes to the natural growth of the population.  If for instance $n(t,x)\geq \beta>0$ throughout the domain, we have that if $u$ is large then the term $\lambda u-n(t,x)|u|^{\rho-1}u<0$ and the population decreases so that we do not expect large values of $u$. 

\par\bigskip

The case where $n(t,x)\equiv n(x)$, so that the equation becomes autonomous,  that is
\begin{equation}\label{eq2intro}\left\{\begin{array}{lr}
u_t-\Delta u = \lambda u - n(x) |u|^{\rho-1}u, & x\in\Omega, \; t>t_0,\\
u=0, & x\in\partial\Omega, \; t>t_0,\\
u(t_0,x)=u_0(x)\geq 0, & x\in \Omega,
\end{array}\right.\end{equation}
has been addressed in the literature when $n(x)\geq \beta>0$ throughout the domain ,see for instance \cite{cantrellcosner, murray}, and also when $n(\cdot)$ vanishes at some region $K_0$ of the domain, that is, 
$$K_0=\{ x\in \bar\Omega: \, n(x)=0\},$$
see \cite{AnibalArrietaPardo, Julian2} and references therein.
Following again the population dynamics interpretation,  the set $K_0$ represents a gifted location where the species is able to grow whereas $\Omega\setminus K_0$ represents a location where the species encounters hardships such as predators, harsh terrain, low resources or intraspecies competition. We refer to $K_0$ as a ``sanctuary'' for the species. Then, by means of the diffusion, part of the growth attained in $K_0$ is migrated outside $K_0$ where the hard conditions threaten the growth of the species. It is natural to see that the size of $K_0$ and the strength of the intrinsic growth rate $\lambda$ will play a key role determining whether the species keeps growing forever or remains bounded.

 Also, in many places in the literature it is considered the case where the function $n(x)$ is regular (at least continuous) and $K_0$ is the closure of a regular open set $\Omega_0$, that is $K_0 = \bar\Omega_0$, see for instance \cite{Julian2} and references therein. With these conditions, the asymptotic behavior of solutions of (\ref{PreEq1}) has been studied. Denoting by $\lambda_1^\Omega$ the first eigenvalue of the Laplace operator with homogeneous Dirichlet boundary conditions on $\Omega$, it is known that if $\lambda < \lambda_1^\Omega$ then $0$ is the unique equilibrium which is globally asymptotically stable and therefore all solutions approach $0$ when $t\rightarrow \infty$. If $\lambda > \lambda_1^\Omega$ then $0$ becomes unstable, and solutions tend to grow away from it. As a matter of fact, if $\lambda_1^\Omega < \lambda < \lambda_1^{\Omega_0}$ then solutions tend to the unique globally asymptotically stable positive equilibrium $u_{eq}$. On the other hand, if $\lambda > \lambda_1^{\Omega_0}$ then solutions have been shown to grow and become unbounded as $t\rightarrow \infty$ only on $K_0$ but they remain bounded outside $K_0$, which will be referred to as ``grow up''. We refer to \cite{Julian2} for a rather complete study of this case. Hence, these results can be summarised as follows:

\begin{equation}\label{introeq2_0}\left\{\begin{array}{ll}
\mbox{If }\lambda < \lambda_1^{\Omega} &\Rightarrow u\rightarrow 0,\\
\mbox{If }\lambda_1^\Omega< \lambda < \lambda_1^{\Omega_0} &\Rightarrow u\rightarrow u_{eq},\\
\mbox{If }\lambda_1^{\Omega_0} < \lambda &\Rightarrow u(x)\rightarrow \infty,\quad x\in \Omega_0.
\end{array}\right.\end{equation}

An extension of (\ref{introeq2_0}) to the case where $K_0$ is only assumed to be compact without any further regularity assumptions was studied in \cite{AnibalArrietaPardo}. In \cite{AnibalArrietaPardo}, the role of $\lambda_1^{\Omega_0}$ is substituted by $\lambda_0(K_0)$, the characteristic value of the set $K_0$, which is defined as follows: consider $\{\Omega_\delta\}_{0<\delta\leq\delta_0}$ a nested decreasing family of smooth open sets satisfyng $K_0=\bigcap_{0<\delta\leq\delta_0} \Omega_\delta$ and define

$$\lambda_0(K_0) := \lim\limits_{\delta \rightarrow 0^+}\lambda_1^{\Omega_\delta}.$$

\pa Notice that $\lambda_0(K_0)$ can be $+\infty$ if $K_0$ is ``small'', for example, when $K_0$ is a point in $\Omega$. If $K_0$ is regular, that is, $K_0 = \bar\Omega_0$ for some regular open set $\Omega_0$, then it is possible to show that  $\lambda_0(K_0) = \lambda_1^{\Omega_0}$. In \cite{AnibalArrietaPardo} authors show that if $\lambda<\lambda_0(K_0)$ then all solutions are bounded and if $\lambda \geq \lambda_0(K_0)$ solutions grow up without bound as time goes to $+\infty$.

\begin{figure}[H]
  \centering
    \includegraphics[width=7cm]{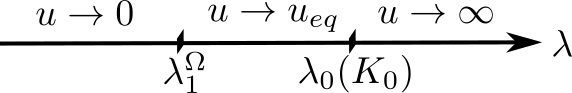}
        \caption{A representation of the cases shown in \eqref{introeq2_0}.}
        \label{figure1}
\end{figure}

We  summarize this results in Figure \ref{figure1}.

\par Therefore, for the autonomous case, the key factor that decides whether the solutions are bounded  or unbounded  is the relation between $\lambda$ and $\lambda_0(K_0)$:  if $\lambda < \lambda_0(K_0)$  solutions remain bounded whereas if $\lambda\geq \lambda_0(K_0)$   solutions become  unbounded as $t\rightarrow \infty$. Moreover, in general terms observe that the smaller $K_0$ is, the larger $\lambda_0(K_0)$ is. Therefore, smaller $K_0$ imply that it is more likely to have $\lambda < \lambda_0(K_0)$ and thus, boundedness. And viceversa, the larger the set $K_0$, the smaller $\lambda_0(K_0)$ is and the more likely that solutions will grow. 
 \bigskip

In the nonautonomous case, that is problem  \eqref{PreEq1},  the situation is not so simple and, as we will see in this work, the relation between $\lambda$ and $\lambda_0(K(t))$ will not be enough to decide the boundedness or unboundedness of the solutions. Note that $\lambda_0(K(t))$ is time dependent and we may not have $\lambda<\lambda_0(K(t))$ for all $t$ or $\lambda>\lambda_0(K(t))$ for all $t$. Moreover, as we will see below, even in the case $\lambda>\lambda_0(K(t))$ for all $t$ we may not have grow--up and if $\lambda < \lambda_0(K(t))$ may not have boundedness.  Therefore, the characterization of boundedness or unboundedness of solutions for the nonautonomous case is not so straightforward  as in the case of autonomous problems. As a matter of fact, many interesting questions appear in the nonautonomous case. Does the geometry of the sets $K(t)$ determine the boundedness character of the solutions?  Is the velocity at which the sets $K(t)$ evolve important to decide the boundedness?  We will see that both the geometry and velocity of the sets $K(t)$ are important factors to take into consideration to decide the boundeness character of the solutions.

\par\bigskip To conclude this introduction, we would like to mention that the results of the present article are a necessary first step to analyze the asymptotic behavior of the solutions of these nonautonomous degenerate logistic equations.  In the cases where boundedness is obtained we should question ourselves about the existence of attractors, their fine structure and behavior.  In the case of unboundedness, we should try to understand how and where the solution becomes unbounded as $t\to +\infty$. Most likely, the set where the solution becomes unbounded (the ``grow-up set'') will be a moving set. These issues will be treated in future publications.  
\par\medskip 
With respect to boundedness, let us mention that there are several interesting works in the literature addressing the general problem of the asymptotic behavior of solutions for nonautonomous problems, analyzing the existence of different kind of asymptotic sets which respond to different concepts of attraction like pullback, forwards, uniform attractors and which describe the behavior of the systems for large times.  Without being exhaustive we would like to mention the monograph \cite{Carvalho2013} and references therein for a nice and rather complete reference in this respect. We also mention the works \cite{anibal2007, tesisAVL} where the authors construct complete trajectores for both autonomous and nonoautonomous parabolic evolution equations, which bound the asymptotic dynamics of the evolution processes. The works  \cite{langa2007, langa2002pullback} are also relevant in this respect.  

Moreover, there are some other interesting works that deal with the periodic case, that is $n(t+T,x)=n(t,x)$ for all $t\in \R$, $x\in \Omega$  and obtain results on existence  of periodic solutions under different hypotheses.  In particular, if  $K(t)$ is static, that is $K(t)\equiv K_0$, existence results of periodic solutions (and therefore, boundedness) were already proven in \cite{Julian3}. Very recently in \cite{Aleja}, authors generalise results in \cite{Julian3} and characterize the existence of periodic solutions without imposing $K(t)$ being static.  In particular, authors prove a very interesting result where they give an if and only if characterization for the existence of periodic solutions in terms of the behavior of certain nontrivial periodic parabolic eigenvalue problems.  This interesting technique is specific for periodic situations and do not apply directly to our most general nonautonomous scenario. See also \cite{DancerHess, Julian4} for some related results on periodic parabolic problems.  

\par\medskip 
With respect to unboundedness, we should try to understand how and where the solution becomes unbounded as $t\to +\infty$. Most likely, the set where the solution becomes unbounded (the ``grow-up set'') will be a moving set but to identify it and study the asympotic properties of these systems requires further analysis.  These issues will be treated in future publications.  

\par\bigskip We describe now the contents of the paper.
\par\medskip
In Section \ref{Section-preliminaries} we include some preliminary results like local and global existence of mild solutions and regularity of solutions. We also recall comparison results with respect to initial data, nonlinearity and domain.  

\par\medskip In Section \ref{Section-data-independence} we include an important result that says that 
the boundedness or unboundedness character of the solution evaluated at a particular point $x\in \Omega$ does not depend on initial time or the initial data $u_0$ (as long as $u_0\geq 0$ and $u_0\not\equiv 0$). Therefore if for one particular nontrivial initial data $u_0$ and initial time $\tau_0$ its solution is asymptotically bounded at $x$, that is $\limsup_{t\to+\infty} u(t,x)<+\infty$, then for any other nontrivial initial data and initial time we also have that the limsup is finite.  This result is based on comparison results, regularity and the parabolic Hopf Lemma.

\par\medskip In Section \ref{Section-general-result} we prove a first result that applies to many situations since we do not impose many restrictions on the configuration of the family of sets $K(t)$. Roughly speaking, the idea is the following: if we assume a nondegeneracy condition of the function $n(t,x)$ (see Assmption (N)) and consider the upper and lower limit of the family of sets $K(t)$, that is $\overline{\bigcup\limits_{t\geq\tau_0}K(t)}$, $\bigcap\limits_{t\geq\tau_0}K(t)$, then if $\lambda<\lambda_0(\overline{\bigcup\limits_{t\geq\tau_0}K(t)})$ we have boundedness of solutions, while if $\lambda >\lambda_0(\bigcap\limits_{t\geq\tau_0}K(t))$ we have unboundeness of solutions, see Proposition \ref{prop-boundedness-general} and Corollary \ref{cor-boundedness-general}.  We also include in this section some relevant examples that allow us to clarify the meaning of these sets.

\par\medskip In Section \ref{secBound} we address situations where the rather general condition of the previous section does not apply and in particular we have $\lambda_0(\overline{\bigcup\limits_{t\geq\tau_0}K(t)})<\lambda<\lambda_0(\bigcap\limits_{t\geq\tau_0}K(t))$. 

A first important case covered in this section says, roughly speaking,  that if $n(t,x)\geq \nu_0>0$  (and therefore $K(t)=\emptyset$) for a sequence of time intervals $I_1$, $I_2$, etc. with the condition that the length of each of them is bounded below by a quantity $\eta>0$ and the distance between two consecutive intervals is no larger than a quantity $\Xi$, then, independently of the value of $\lambda$, all solutions are bounded, see Proposition \ref{PropInter}.   That is, the effective presence of the logistic term in the entire domain and at the intervals of time $I_1$, $I_2$, etc.  which do not degenerate and they are not too far away from the next, guarantees the boundedness of solutions.

\begin{figure}[H]
  \centering
    \includegraphics[width=8cm]{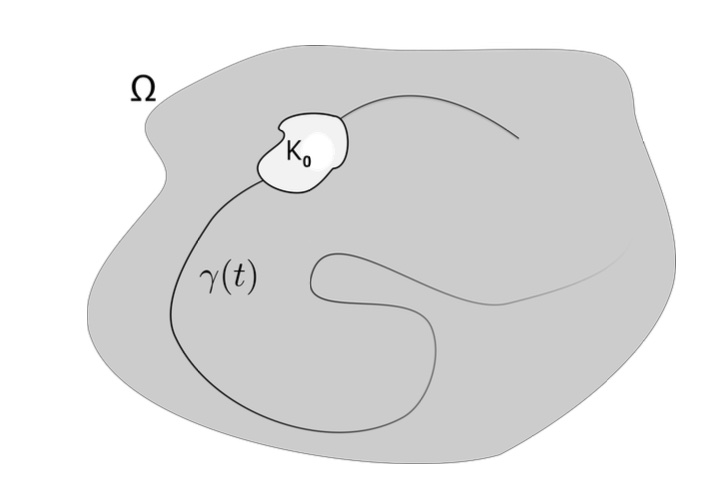}
    \caption{Domain $\Omega$ and a continuously moving $K(t)$. }
    \label{FigPeriodicTranslating}
\end{figure}

Next, we show that, roughly speaking and oversimplifying the conditions,  if $\lambda<\lambda_0(K(t))$ for all $t$ large enough  then we have boundedness of solutions. The precise statement is contained in Theorem \ref{TeoAcot1}.  A very important example where this theorem applies is contained in Remark \ref{RecondAcotLambdaSmall} where a fix set $K_0$ is considered and we assume that $K(t)$ consists in moving (translation by a curve $\gamma(t)$ and a rotation) the set $K_0$ inside $\Omega$ see Figure \ref{FigPeriodicTranslating}. 

In this section we also consider the case where the set $K(t)$ allows jumps. A good example of this situation is when we consider the existence of a sequence $t_0<t_1<t_2<\ldots$, such that $K(t)=K_0$ for $t\in (t_{2k},t_{2k+1})$ and $K(t)=K_1$ for $t\in (t_{2k+1},t_{2k+2})$.  Hence, the set $K(t)$ jumps from $K_0$ to $K_1$ and back to $K_0$ at certain instants of time. We show that if $K_0\cap K_1=\emptyset$ then  independently of the value of $\lambda$ all solutions are bounded.  

Finally, in Section \ref{secUnbound} we analyze conditions under which solutions become unbounded.   In view of the results of Theorem \ref{TeoAcot1} where boundedness is established for, roughly speaking, $\lambda<\lambda_0(K(t))$ for all $t$,  it seems natural to ask whether for $\lambda>\lambda_0(K(t))$ for all $t$,  we have unboundedness of solutions. This is actually true if the equation is autonomous, see \cite{AnibalArrietaPardo} for instance. But notice that Theorem \ref{TeoAcot} and the example where $K(t)$ jumps between two sets $K_0$ and $K_1$ with $K_0\cap K_1=\emptyset$ and $\lambda$ arbitrary and in particular we could consider $\lambda>\lambda_0(K_0), \lambda_0(K_1)$, means that this is not true in a general setting.  As a matter of fact, we will be able to show that if, roughly speaking,  $K(t)$ moves ``slowly'' and $\lambda>\lambda_0(K(t))$ for all $t$, then we have unboundedness of solutions.  The precise statement of this result is Theorem \ref{TeoUnbounded}.  We also analyze a situation where the condition $\lambda>\lambda_0(K(t))$ for all $t$ can be weakened, see Proposition \ref{PropUnbound}.

\par\bigskip\noindent {\bf Acknowledgements.}  José M. Arrieta is partially supported by grants PID2019-103860GB-I00, PID2022-137074NB-I00  and CEX2019-000904-S ``Severo Ochoa Programme for Centres of Excellence in R\&D'' the three of them from MICINN, Spain. Also by ``Grupo de Investigación 920894 - CADEDIF'', UCM, Spain.   Marcos Molina-Rodríguez is partially supported by the predoctoral appointment  BES-2013-066013, grant PID2019-103860GB-I00 both by MICINN, Spain and by ``Grupo de Investigación 920894 - CADEDIF'', UCM, Spain.  Lucas A. Santos is partially supported by the Coordenação de Aperfeiçoamento de Pessoal de Nível Superior – Brasil (CAPES) – Finance Code 001”

\section{Preliminaries}
\label{Section-preliminaries}

Equation (\ref{PreEq1}) is a particular case of a more general class of parabolic nonautonomous equations:

\begin{equation}\label{PreEq2}\left\{\begin{array}{lr}
u_t- \Delta u = f(t,x,u), & x\in\Omega, \; t>t_0,\\
u=0, & x\in\partial \Omega, \; t>t_0,\\
u(t_0,x)=u_0(x), & x\in \Omega,
\end{array}\right.\end{equation}

\pa where $f:\R\times\Omega\times\R\rightarrow \R$ satisfies

\begin{equation}\label{condfExUn}\begin{array}{l}
\left| f(t,x,u) - f(t,x,v)\right| \leq C \left| u - v \right| \left( \left|u\right|^{\rho-1} + \left|v\right|^{\rho-1} + 1\right),  \; \mbox{ for } u,v\in \R, \,\, t\geq t_0,
\end{array}\end{equation}

\pa for certain constant $C$ uniformly for $t\geq t_0$ on bounded sets and  $x\in \Omega$.  Observe that if $n\in L^\infty (\R\times\Omega)$ we easily show that $f(t,x,u)=\lambda u-n(t,x)|u|^{\rho-1}u$ satisfies \eqref{condfExUn} with $C=C(\lambda,\rho,\|n\|_{L^\infty})$.

\medskip

\subsection{Local existence}\label{sec11}

For problems of the type (\ref{PreEq2}) with $u_0\in L^q(\Omega)$ and for $1<q<\infty$ large enough we have existence and uniqueness of local mild solutions, that is, solutions of the variations of constants formula associated to \eqref{condfExUn}. As a matter of fact, following the general techniques to obtain well posedness by \cite{henry} and in particular the results from \cite{Arr2000} we can show that if $q>\D\max\left\{\frac{N}{2}\rho, \rho\right\}$ then problem (\ref{PreEq2}) has a unique $u\in C([t_0,\tau_0),L^q(\Omega))$, local mild solution for some $\tau_0>t_0$. Moreover, for any $0<\eta<1$  $u\in C\left((t_0,\tau_0) ;  L^\infty(\Omega)\right) \cap C\left((t_0,\tau_0); C^{1,\eta}(\bar\Omega)\right)$.

We refer to \cite{TesisMarcos} for more details on this.

\subsection{Positivity, Comparison and Global Existence of Solutions}\label{sec12}

Comparison results for mild solutions as studied in \cite{ArrietaCarvalhoRodriguez} can be applied to (\ref{PreEq2}) and in particular to  (\ref{PreEq1}) with $f(t,x,u)=\lambda u - n(t,x)|u|^{\rho-1}u$.  Let us denote by $u(t,x,t_0,u_0,f,\Omega)$ the unique solution of (\ref{PreEq2}) where  we make explicit the dependence of the solution in the initial time $t_0$ and initial data $u_0$, the nonlinearity $f$ and the domain $\Omega$.  When we want to stress the dependence of the solution on just the initial data and nonlinearity, we may write $u(t,x, u_0,f)$. 

With the techniques and results from \cite{ArrietaCarvalhoRodriguez} it can be easily shown the following 

\lem\label{CorComp} Let $f_i(t,x,u)=\lambda u - n_i(t,x) |u|^{\rho-1} u$ with $n_i\in L^\infty(\R\times  \Omega)$, $n_1(t,x)\geq n_2(t,x)\geq 0$  a.e. $t\in (t_0,+\infty)$, a.e. $x\in \bar\Omega$. Then, for as long as solutions exist, we have that:

\begin{equation}\label{EqCorComp}
0\leq \left|u(t,x,u_0,f_1)\right| \leq u(t,x,|u_0|,f_1) \leq u(t,x,|u_0|,f_2), \quad x\in \bar\Omega.
%= e^{\lambda (t-t_0)}e^{\Delta (t-t_0)} |u_0| ,
\end{equation}
Moreover, if $\tilde\Omega\subset \Omega$ then for any initial condition $u_0\geq 0$ defined in $\Omega$ if we denote by $\tilde u_0= {u_0}_{|\tilde\Omega}$, we have
$$0\leq u(t,x,\tilde u_0, f,\tilde \Omega)\leq u(t,x, u_0,f, \Omega),\quad x\in \tilde \Omega,$$
as long as both solutions exist.

\rm

\par\bigskip
 From the previous Lemma \ref{CorComp} and taking $n_2(t,x)\equiv 0$, we have 
\begin{equation}\label{PreDesuU}
0\leq |u(t,x)| \leq U(t,x),
\end{equation}
as long as the solutions exist and a.e. $x\in \Omega$,  where we have denoted by $U$ the solution of the linear problem
\begin{equation}\label{eqUU}
\left\{ \begin{array}{lr}
U_t - \Delta U = \lambda U , & x\in \Omega,\; t\in(t_0,\infty),\\
U=0, & x\in\partial\Omega, \; t\in(t_0,\infty),\\
U(t_0)=|u_0|, & x\in \Omega. 
\end{array}\right.
\end{equation}

\pa This solution can be expressed also as the linear semigroup $U(t,\cdot)=e^{\lambda (t-t_0)} e^{\Delta(t-t_0)}|u_0|$, which is globally defined and uniformly bounded in compact sets of $t>t_0$. 

\cor[Global Existence]\label{ExistenciaGlobal} With the assumptions above, we have that $u(t,x,t_0,u_0)$ is globally defined and bounded in compact sets of $t\in (t_0,+\infty)$. In particular, for $0\leq u_0\in L^q(\Omega)$ with  $q>\D\max\left\{\frac{N}{2}\rho, \rho\right\}$, $u$ verifies
\begin{equation*}
0\leq u(t,x,t_0,u_0) \leq e^{\lambda t} e^{\Delta (t-t_0)} u_0(x),
\end{equation*}

\pa which implies that

\begin{equation}\label{cotaU}\left\{ \begin{array}{l}
\|u(t,\cdot,u_0)\|_{L^\infty(\Omega)} \leq M e^{(\lambda -\lambda_1^\Omega)(t-t_0)}\|u_0\|_{L^\infty(\Omega)},\quad t>t_0,\\ \\
\|u(t,\cdot,u_0)\|_{L^\infty(\Omega)} 
\leq M e^{(\lambda -\lambda_1^\Omega)(t-t_0)}((t-t_0)^{-\frac{N}{2q}}+1) \|u_0\|_{L^q(\Omega)},\quad t>t_0, %1\leq q\leq\infty
\end{array}\right.\end{equation}
\pa where the constant $M$ is defined so that 

\begin{equation}\label{estU}
\left\{
\begin{array}{l}
\|e^{\Delta t}u_0\|_{L^\infty(\Omega)} \leq M 
e^{-\lambda_1^\Omega t}\|u_0\|_{L^\infty(\Omega)},\quad t>0, 
\\ \\
\|e^{\Delta t}u_0\|_{L^\infty(\Omega)} \leq M 
e^{-\lambda_1^\Omega t}(t^{-\frac{N}{2q}}+1) \|u_0\|_{L^q(\Omega)},\quad t>0.
\end{array}
\right.
\end{equation}
\rm

\re Observe that the constant $M$ above depends only on  $\lambda_1^\Omega$ and $|\Omega|$, see \cite{Cazenave-Haraux}.
\rm

\section{Initial Data Independence}
\label{Section-data-independence}

   We prove now that the boundedness or unboundedness character of the solutions is independent of initial data and initial time. This was shown for the autonomous case in \cite{AnibalArrietaPardo}.  This means that if for some particular initial data $u_0$ and initial time $t_0$ the solution is bounded at certain point $x\in \Omega$, then for any other initial data and initial time, the solution will also be bounded at this point. 
   
   \par\medskip  Let us start with the following auxiliary lemma.
   
\lem\label{LemDat} Let $u(t,x,u_0)$ be the solution of (\ref{PreEq1}) and let $\lambda\in \R$, $u_0\geq 0$  and $\alpha>0$. We have:
\par\noindent i) If  $\alpha>1$, then $\alpha u(t,x,u_0)\geq u(t,x,\alpha u_0).$
\par\noindent ii) If $\alpha<1$, then $\alpha u(t,x,u_0)\leq u(t,x, \alpha u_0).$
\proof The function $v(t,x)=\alpha u(t,x, u_0)$ satisfies  the problem

\begin{equation}
\left\{\begin{array}{lr}
v_{t}-\Delta v=\lambda v - \frac{n(t,x)}{\alpha^{\rho-1}} v^\rho,& x\in \Omega, t>0,\\
v(t,x)=0, & x\in \partial \Omega, t>0,\\
v(0,x)=\alpha u_0,& x\in \Omega.
\end{array}
\right.\end{equation}

By comparison results, see Lemma \ref{CorComp}, since $\frac{n(t,x)}{\alpha^{\rho -1}}\leq n(t,x)$ if $\alpha>1$ and $\frac{n(t,x)}{\alpha^{\rho -1}}\geq n(t,x)$ if $\alpha<1$, we have

\par\noindent i) if $\alpha>1$, then $v$ is a supersolution of $u(t,x,\alpha u_0)$, so $\alpha u(t,x,u_0) \geq u(t,x,\alpha u_0)$.

\par\noindent ii) if $\alpha<1$, then $v$ is a subsolution of $u(t,x,\alpha u_0)$, so $\alpha u(t,x,u_0)\leq u(t,x,\alpha u_0)$.

\par This proves the result.\eproof

\par\medskip
Before stating and proving the main result of this section, let us state a Hopf parabolic Lemma for our equation \eqref{PreEq1}  and a technical result. We we will prove both of them below, after the main result.

\lem\label{Hopf} (Hopf parabolic Lemma) Let $\Omega\subset \R^N$ satisfy the interior sphere condition, $n\in L^\infty(\R\times \Omega)$ and let $u(t,x,t_0,u_0)$ be the solution of (\ref{PreEq1}) with $u_0\in L^\infty(\Omega)$,  $u_0\geq 0$ and $u_0\not\equiv 0$. Then we have  that

\begin{equation*}\label{HopfU}
\frac{\partial u}{\partial \vec n} (t,x,u_0) < 0, \; \; \;  \forall t>t_0, x\in \partial \Omega,
\end{equation*}

\pa where $\vec n$ is the outward normal vector of $\partial \Omega$.
\rm 
\par  Recall that a domain $\Omega$ satisfies the interior sphere condition if for each $x\in \partial\Omega$ there exists $z\in\Omega$ such that if $a=|x-z|$ then $B(z,a)\subset \Omega$ and $\bar B(z,a) \cap \bar\Omega=\{ x\}$.

\lem \label{auxiliar-Hopf}Let $\Omega$ be as above and let $\varphi,\phi\in C^1(\bar\Omega)$ with $\varphi(x),\phi(x)>0$ for $x\in\Omega$ and $\varphi(x)=\phi(x)=0$  for $x\in\partial\Omega$. Then if $\frac{\partial\varphi}{\partial n}(x)<0$ for each $x\in \partial \Omega$, there exists 
$\alpha<1$ such that $$\alpha\phi(x)\leq \varphi(x), \hbox{  for all }x\in \bar\Omega.$$
\rm 

\par\medskip  Now, we can prove,

\prop\label{IndepDatIni} The solution of (\ref{PreEq1}) is bounded or unbounded independently of initial data and initial time.

\proof Let $u_0,v_0$ be two non--trivial, non--negative initial data, let $\tau_0\geq t_0$  and let us denote by $u(t,x;t_0,u_0)$ be the solution of (\ref{PreEq1}) starting at $t_0$ with initial condition $u_0$. 

First, we want to show that $u(t,\cdot; \tau_0,u_0)$ is bounded for $t>\tau_0$  iff $u(t,\cdot;\tau_0,v_0)$ is bounded for $t>\tau_0$ . We let the process evolve a small amount of time $\delta>0$ and by the regularity stated in Subsection \ref{sec11},   $u(\tau_0+\delta,\cdot;\tau_0, u_0),u(\tau_0+\delta,\cdot;\tau_0, v_0)\in C^1(\bar\Omega)$.  Moreover, these solutions vanish in $\partial\Omega$, both are strictly positive inside $\Omega$ and by the parabolic Hopf Lemma (see Lemma \ref{Hopf} )  they also satisfy 
$$\frac{\partial u(\tau_0+\delta,\cdot;\tau_0, u_0)}{\partial \vec n}<0,\qquad \frac{\partial u(\tau_0+\delta,\cdot;\tau_0, v_0)}{\partial \vec n}<0.$$

Applying now Lemma \ref{auxiliar-Hopf} to these two functions, $x\to u(\tau_0+\delta,x;\tau_0, u_0)$ and $x\to u(\tau_0+\delta,x;\tau_0, v_0)$ both playing the role of $\varphi$ and $\phi$, we have the existence of $\alpha,\tilde \alpha<1$ such that 
$$\alpha u(\tau_0+\delta,\cdot;\tau_0,u_0)\leq u(\tau_0+\delta,\cdot;\tau_0,v_0),  \quad  \tilde\alpha u(\tau_0+\delta,\cdot;\tau_0,v_0)\leq u(\tau_0+\delta,\cdot;\tau_0,u_0).$$

 Therefore, denoting  $\beta=1/\tilde\alpha$, we get that 

$$\alpha u(\tau_0+\delta,\cdot;\tau_0,u_0)\leq u(\tau_0+\delta,\cdot;\tau_0,v_0) \leq \beta u(\tau_0+\delta,\cdot;\tau_0,u_0).$$

\pa Hence, taking into account that $u(t,\cdot;\tau_0,v_0)=u(t,\cdot;\tau_0+\delta,u(\tau_0+\delta,\cdot;\tau_0,v_0))$, for $t>\tau_0+\delta$ we have

\begin{equation*}
 u(t,\cdot;\tau_0+\delta,\alpha u(\tau_0+\delta,\cdot;\tau_0,u_0)) \leq u(t,\cdot;\tau_0,v_0) \leq  u(t,\cdot;\tau_0+\delta,\beta u(\tau_0+\delta,\cdot;\tau_0,u_0)).
\end{equation*}

\pa But from Lemma \ref{LemDat} for $\alpha<1<\beta$ and $t> t_0 + \delta$
\begin{equation*}\begin{array}{c}
\alpha  u(t,\cdot;\tau_0,u_0) = \alpha  u(t,\cdot;\tau_0+\delta,u(\tau_0+\delta,\cdot;\tau_0,u_0)) \leq  u(t,\cdot;\tau_0+\delta,\alpha u(\tau_0+\delta,\cdot;\tau_0,u_0)), \\
\; \\
 u(t,\cdot;\tau_0+\delta,\beta u(\tau_0+\delta,\cdot;\tau_0,u_0)) \leq \beta  u(t,\cdot;\tau_0+\delta,u(\tau_0+\delta,\cdot;\tau_0,u_0))= \beta  u(t,\cdot;\tau_0,u_0).
\end{array}\end{equation*}

\pa Hence, we arrive at

\begin{equation*}
\alpha u(t,\cdot;\tau_0,u_0)\leq u(t,\cdot;\tau_0,v_0) \leq \beta u(t,\cdot;\tau_0,u_0), \quad t>\tau_0+\delta.
\end{equation*}

\pa Thus, one solution is bounded if and only if the other is bounded. 

\pa Now let us take two different initial times $\tau_1,\tau_2\geq t_0$ and two initial data $u_0, v_0$ as before. Without any loss of generality, we may assume $\tau_1<\tau_2$. Let the process evolve from $\tau_1$ to $\tau_2$ with initial data $u_0$ and denote by $\tilde u_0=u(\tau_2,\cdot;\tau_1,u_0)$. From the proof above, we know that
$$u(t,\cdot;\tau_2,\tilde u_0)=u(t,\cdot;\tau_1,u_0) \quad \mbox{is bounded for } t\geq\tau_2 \quad\Leftrightarrow \quad u(t,\cdot;\tau_2,v_0) \mbox{ is bounded}.$$
This shows the result. 
\eproof

\medskip  For the sake of completeness we include a short proof of Lemma \ref{Hopf}. 

\proof (Hopf parabolic Lemma)  Denoting by $U$ the solution of the linear problem (\ref{eqUU}), by Lemma \ref{CorComp} we have that for $t\geq t_0$, $U(t) \in L^\infty(\Omega)$ and $\|U(t)\|_\infty \leq e^{\lambda (t-t_0)} \|u_0\|_\infty$. Moreover,

\begin{equation*}
0\leq u(t,x; u_0) \leq U(t,x;u_0), \; \mbox{ for } t>t_0, x\in\bar\Omega.
\end{equation*}

\pa Then,

\begin{equation*}
\lambda u - n(t,x) |u|^{\rho-1} u \geq \lambda u - \|n\|_\infty \left(e^{\lambda (t-t_0)} \|u_0\|_\infty\right)^{\rho-1} u \geq \lambda u - N_0 u, \quad t\in[t_0,t_1],
\end{equation*}

\pa where 

\begin{equation*}
N_0=\D\sup_{t\in [t_0,t_1]}\|n\|_\infty \left(\|u_0\|_\infty e^{\lambda (t-t_0)}\right)^{\rho-1} =\|n\|_\infty \left(\|u_0\|_\infty e^{|\lambda| (t_1-t_0)} \right)^{\rho-1} .
\end{equation*}

\pa This implies that the solution $v$ of
\begin{equation}\label{EqVhopf}\left\{\begin{array}{lr}
v_t-\Delta v-(\lambda-N_0)v=0 & t\in (t_0,t_1),\; x\in \Omega,\\
v=0, & t\in(t_0,t_1),\;  x\in\partial \Omega,\\
v(t_0)=u_0 & x\in \Omega, 
\end{array}\right.\end{equation}

\pa satisfies, using Lemma \ref{CorComp}, $u(t,x)\geq v(t,x)$. Moreover, multiplying $e^{-(\lambda-N_0)(t-t_0)}$ to both sides of equation for $v$, we arrive at $v(t,x)=e^{(\lambda-N_0)(t-t_0)}w(t,x)$ where

\begin{equation}\label{EqWhopf}\left\{\begin{array}{lr}
w_t-\Delta w=0 & (t,x)\in (t_0,t_1)\times \Omega,\\
w=0, & x\in\partial \Omega,\\
w(t_0)=u_0 & x\in \Omega.
\end{array}\right.\end{equation}

\pa Therefore,

\begin{equation}\label{DesUW}
u(t,x)\geq v(t,x) = e^{(\lambda - N_0)t} w(t,x), \; \mbox{ for all } (t,x) \in (t_0,t_1)\times \Omega.
\end{equation}

Since $u(t,x)=v(t,x)=w(t,x)=0$ at the boundary and by the regularity of $u(t,\cdot)$ we have

\begin{equation}\label{desigualdadesHopf}
\frac{\partial u}{\partial \vec n}(t,x)\leq \frac{\partial v}{\partial \vec n}(t,x)=e^{(\lambda - N_0)t}\frac{\partial w}{\partial \vec n}(t,x), \; \mbox{ for all } (t,x) \in (t_0,t_1)\times \partial\Omega.\end{equation}

  Applying the  Hopf parabolic Lemma to $w$ which is a classical solution, see  \cite[Chapter 3 Theorem 7]{protter}, we get $\frac{\partial w}{\partial \vec n} (t,x,u_0) < 0$  for all $t>t_0$, $x\in \partial \Omega$.  This statement together with \eqref{desigualdadesHopf} imply the result. \eproof

\proof (Lemma \ref{auxiliar-Hopf}) Let us denote by $A>0$ a constant such that $|\nabla \varphi(x)|,|\nabla\phi(x)|, \varphi(x),\phi(x)\leq A$ for all $x\in\bar\Omega$. Let us argue by contradiction. If the result is not true, then there exists $\alpha_n\to 0$ and $x_n\in \bar\Omega$ such that $\alpha_n\phi(x_n)>\varphi(x_n)$.  
This implies that $x_n\in \Omega$,  $\varphi(x_n)\leq A\alpha_n\to 0$ and therefore $d(x_n,\partial\Omega)\to 0$. Taking a subsequence if necessary we obtain the existence of $x^*\in\partial\Omega$ such that $x_n\to x^*$.   Since $\frac{\partial\varphi}{\partial n}(x)<0$, by continuity of $x\to \nabla\varphi(x)$ in $\bar\Omega$ and $x\to \vec n(x)$ in $\partial\Omega$, we get
that there exists $\delta, a>0$ such that 
\begin{equation}\label{auxiliar01}
\nabla\varphi(x)\cdot \vec n(y)\leq -a, \quad \forall x\in B(x*,\delta)\cap \bar\Omega,\quad \forall y\in B(x*,\delta)\cap \partial\Omega.
\end{equation}

Since $x_n\to x^*$ we have that there exists $n_0\in\N$ such that $|x_n-x^*|<\delta/2$ for $n\geq n_0$. If we take now $\rho_n$ such that $B(x_n,\rho_n)\subset \Omega$ and $\bar B(x_n,\rho_n)\cap\partial\Omega\ni z_n^*$ then we have $\rho_n<\delta/2$ and $|z_n^*-x^*|\leq |z_n^*-x_n|+|x_n-x^*|<\delta/2+\delta/2=\delta$.  

Using \eqref{auxiliar01} and the mean value theorem, we have:

$$\frac{\varphi(x_n)-\varphi(z_n^*)}{|x_n-z_n^*|}=\nabla \varphi(\xi_n)\cdot \frac{x_n-z_n^*}{|x_n-z_n^*|}=-\nabla\varphi(\xi_n)\vec n(z_n^*)\geq a>0, \quad \forall n\geq n_0,$$
where $\xi_n$ lies in the segment joining $x_n$ and $z_n^*$. Observe that we have used that the unit exterior normal vector at $z_n^*$ is given $-\frac{x_n-z_n^*}{|x_n-z_n^*|}$.

On the other hand, 
$$\frac{\varphi(x_n)-\varphi(z_n^*)}{|x_n-z_n^*|}=\frac{\varphi(x_n)}{|x_n-z_n^*|}\leq \alpha_n\frac{\phi(x_n)}{|x_n-z_n^*|}=\alpha_n\frac{\phi(x_n)-\phi(z_n^*)}{|x_n-z_n^*|}\leq \alpha_nA\to 0,$$
which is in contradiction with the inequality just obtained above.\eproof

\section{A general result}\label{Section-general-result}

\medskip
 
We start giving in this section a rather general result, which can be applied in many situations with not many geometric restrictions on the sets $K(t)$, where boundedness or unboundedness is proved. These results are obtained  via comparison techniques with certain autonomous reference problems involving,  roughly speaking, the upper and lower limit of the sets $K(t)$. 
\par\medskip 
 \par Moreover,  we will often assume the following condition:

\begin{assf*}[\bf N] There exists a strictly increasing continuous function $\nu:\R^+\rightarrow \R^+$ with $\nu(0)=0$ and
\begin{equation*}\label{nt}
n(t,x)\geq \nu(d(x,K(t))),\; \mbox{ for all } x\in \overline{\Omega}, t\in \R.
\end{equation*} 

\pa Notice that $K(t)$ may be empty for some value of $t$. For the values of $t$ when $K(t)$ is empty, we will assume that there exists $\nu_0>0$ such that $n(t,x)\geq \nu_0$ for all $x\in\overline{\Omega}$.

\end{assf*}

Let us define the following functions and sets which will be important in the analysis below.  For a fixed $\tau_0\geq t_0$, 
\begin{equation}\label{ntau0}
\begin{array}{l}
\underline{n}_{\tau_0}(x):=\inf\limits_{t\geq\tau_0} n(t,x), \; \mbox{ for }x\in\overline{\Omega}, \qquad \uK_{\tau_0}:=\overline{\{ x\in \overline{\Omega} : \underline{n}_{\tau_0}(x)=0\}}, %\quad \tau_0\geq t_0   
\\ \\ 
\overline{n}_{\tau_0}(x):=\sup\limits_{t\geq\tau_0} n(t,x), \; \mbox{ for }x\in\overline{\Omega},
\qquad \oK_{\tau_0}:=\{ x\in \overline{\Omega} : \overline{n}_{\tau_0}(x)=0\}. %\quad \tau_0\geq t_0 
\end{array}
\end{equation}

We  can show, 
\lem\label{KsetDef} With the notations above, $\uK_{\tau_0}$ is decreasing in $\tau_0$, $\oK_{\tau_0}$ is increasing in $\tau_0$ and 
\begin{equation*}
\overline{\bigcup\limits_{t\geq\tau_0}K(t)}\subset \uK_{\tau_0},
\qquad \bigcap\limits_{t\geq\tau_0}K(t)=\oK_{\tau_0}.
\end{equation*}

\pa  Moreover if {\bf(N)} holds, then $\overline{\bigcup\limits_{t\geq\tau_0}K(t)}= \uK_{\tau_0}$.
\proof That $\uK_{\tau_0}$ is decreasing in $\tau_0$, $\oK_{\tau_0}$ is increasing in $\tau_0$  follows directly from the definition in \eqref{ntau0}.  

Let us see all inclusions. Let $x\in \bigcup\limits_{t\geq\tau_0} K(t)$, then  $x\in K(t)$ for some $t\geq\tau_0$, which implies $n(t,x)=0$ and therefore $\underline{n}_{\tau_0}(x)=0$. Hence, 

\begin{equation*}
\bigcup\limits_{t\geq\tau_0} K(t)\subset \{ x\in \Omega : \underline{n}_{\tau_0}(x)=0\} \Rightarrow \overline{\bigcup\limits_{t\geq\tau_0} K(t)}\subset \overline{\{ x\in \Omega : \underline{n}_{\tau_0}(x)=0\}}=\uK_{\tau_0}.
\end{equation*}

Now, let us take $x\in \oK_{\tau_0}$, then $\sup\limits_{t\geq\tau_0} n(t,x)=0$. Therefore

\begin{equation*}
0\leq n(t,x) \leq \sup\limits_{t\geq\tau_0} n(t,x)=0, \;\; \forall t\geq\tau_0\;\; \Rightarrow n(t,x)=0,\;\; \forall t\geq\tau_0.
\end{equation*}

\pa Hence, $x\in K(t)$ for every $t\geq\tau_0$ and thus, $x\in\bigcap\limits_{t\geq\tau_0} K(t)$ and therefore, $\oK_{\tau_0}\subset\bigcap\limits_{t\geq\tau_0}K(t)$.

Now, let $x\in \bigcap\limits_{t\geq\tau_0} K(t)$, then $x\in K(t)$ for every $t\geq\tau_0$. Therefore, $\overline{n}_{\tau_0}(x)=0$ which implies that $x\in \oK_{\tau_0}$ and, thus, $\bigcap\limits_{t\geq\tau_0}K(t)\subset\oK_{\tau_0}$.

Finally, let us assume {\bf(N)} holds and $x\in \{ x\in \Omega : \underline{n}_{\tau_0}(x)=0\}$ for some $\tau_0\geq t_0$. Then, $\underline{n}_{\tau_0}(x)=\inf\limits_{t\geq\tau_0} n(t,x)=0$. Therefore, either 

\begin{itemize}

\item $\exists t\geq \tau_0$ with $n(t,x)=0$. This implies $x\in K(t)$ and therefore $x\in \bigcup\limits_{t\geq\tau_0} K(t)$ or,

\item $n(t,x)>0$ for every $t\geq\tau_0$ and there exists a sequence $(t_i)_{i\in\N}$ such that $\lim\limits_{i\rightarrow\infty} n(t_i,x)=0$, which implies that
\begin{equation*}
0\leq \lim\limits_{i\rightarrow\infty} \nu(d(x,K(t_i))) \leq \lim\limits_{i\rightarrow\infty} n(t_i,x)=0.
\end{equation*}

\pa Hence, $\lim\limits_{i\rightarrow \infty} d(x,K(t_i)) = 0$ and, therefore,
 $x\in  \overline{\bigcup\limits_{t\geq\tau_0} K(t)}$.
\end{itemize}

\pa Thus, $\uK_{\tau_0}\subset \overline{\bigcup\limits_{t\geq\tau_0} K(t)}$. This concludes the proof of the result.  \eproof

\re\label{Rent} Assumption  {\bf(N)} is avoiding that the function $n(t,x)$ degenerates outside $K(t)$, that is, if $x(t)\in \Omega$ is at a positive distance from $K(t)$, that is, $d(x(t),K(t))\geq d_0>0$, then $n(t,x(t))\geq \nu(d_0)>0$ for all $t$.   This is important  for Lemma \ref{KsetDef} and also for  many  results in this work. Notice that we can construct the function 
\begin{equation*}
n(t,x)=\left\{\begin{array}{lr}
0,& t\geq t_0,x\in K_0, \\
\frac{1}{t-t_0+1}, & t\geq t_0, \; x\in\overline{\Omega}\setminus K_0,
\end{array}\right.
\end{equation*}

\pa where $K(t)\equiv K_0$ compact but $\underline{n}_{\tau_0}(x)=0$ for all $\tau_0\geq t_0$ and all $x\in\overline{\Omega}$.  Therefore $\uK_{\tau_0}=\bar\Omega$ but $\overline{\bigcup\limits_{t\geq \tau_0} K(t)}=K_0 \subsetneq \bar\Omega$, and therefore the last inequality in Lemma \ref{KsetDef} does not hold true.\rm

\medskip

\medskip

\rm We can show now

\prop \label{prop-boundedness-general} Let us consider problem (\ref{PreEq1}) and let us assume {\bf (N)} holds.  We have, 
\begin{itemize} 
\item If $\lambda<\lambda_0(\uK_{\tau_0})$ for some $\tau_0\geq t_0$, then all solutions of (\ref{PreEq1}) are bounded.
\item If  $\lambda\geq\lambda_0(\oK_{\tau_0})$ for some $\tau_0\geq t_0$, then all solutions of (\ref{PreEq1}) are unbounded.
\end{itemize}

\proof  Let ${\tau_0}\geq t_0$ and define 
\begin{equation*}\underline{n}_{\tau_0}^*:=\left\{\begin{array}{lr}
0, & x\in \uK_{\tau_0},\\
\underline{n}_{\tau_0}(x), &x\in\overline{\Omega}\setminus\uK_{\tau_0}.
\end{array}\right.\end{equation*}

Let $\underline{v}$ and $\overline{v}$ be the respective solutions of the autonomous problems

\begin{equation*}
\left\{\begin{array}{ll}
\underline{v}_t-\Delta \underline{v} = \lambda \underline{v} - \underline{n}_{\tau_0}^* |\underline{v}|^{\rho-1}\underline{v}, & x\in\Omega, \; t>\tau_0,\\
\underline{v}=0, & x\in\partial\Omega, \; t>\tau_0,\\
\underline{v}(\tau_0,x)=u_0(x)\geq 0, &  x\in\Omega,
\end{array}\right.
\quad 
\left\{\begin{array}{ll}
\overline{v}_t-\Delta \overline{v} = \lambda \overline{v} - \underline{n}_{\tau_0} |\overline{v}|^{\rho-1}\overline{v}, & x\in\Omega, \; t>\tau_0,\\
\overline{v}=0, & x\in\partial\Omega, \; t>\tau_0,\\
\overline{v}(\tau_0,x)=u_0(x)\geq 0, &  x\in\Omega.
\end{array}\right.
\end{equation*}

\pa Then, since 

\begin{equation*}
\overline{n}_{\tau_0}(x)\geq n(t,x)\geq \underline{n}_{\tau_0}(x)\geq \underline{n}_{\tau_0}^*(x), \;  t\geq\tau_0, \; x\in \overline{\Omega},
\end{equation*}

\pa by comparison, see Corollary \ref{CorComp},  we have that

\begin{equation*}
\overline{v}(t,x) \leq u(t,x) \leq \underline{v}(t,x), \; t\geq\tau_0, \mbox{ a.e. } x\in \Omega.
\end{equation*}

\pa Now, 
\begin{itemize}
\item If $\lambda<\lambda_0(\uK_{\tau_0})$ then $\underline{v}$ is bounded due to \cite{AnibalArrietaPardo}. Thus $u$ is a bounded solution.
\item If $\lambda\geq \lambda_0(\oK_{\tau_0})$, then $\overline{v}$  is unbounded due to \cite{AnibalArrietaPardo}. Thus $u$ is an unbounded solution.\eproof
\end{itemize} 

\medskip

If we define $\lambda_0^{-}=\lim_{\tau_0\to+\infty}\lambda_0(K^{\sup}_{\tau_0})\in \R^+\cup\{+\infty\}$ and $\lambda_0^{+}=\lim_{\tau_0\to+\infty}\lambda_0(K^{\inf}_{\tau_0})\in \R^+\cup\{+\infty\}$, which both limits exist since the sequences  $\lambda_0(K^{\sup}_{\tau_0})$ and $\lambda_0(K^{\inf}_{\tau_0})$ are monotone decreasing and increasing respectively, we can obtain

\cor\label{cor-boundedness-general} In the assumptions of Proposition \ref{prop-boundedness-general}, we have, 
\begin{itemize} 
\item If $\lambda<\lambda_0^{-}$ all solutions of (\ref{PreEq1}) are bounded.
\item If  $\lambda>\lambda_0^{+}$ all solutions of (\ref{PreEq1}) are unbounded.
\end{itemize}
\rm

\re  It is natural to define  $K^{\sup}=\cap_{\tau_0\geq t_0}K^{\sup}_{\tau_0}$ and $K^{\inf}=\overline{\cup_{\tau_0>t_0}K^{\inf}_{\tau_0}}$ and these two sets are like the upper and lower limits of the one parameter family of sets $K(t)$.  We would expect that $\lambda_0^-=\lambda_0(K^{\sup})$ and $\lambda_0^+=\lambda_0(K^{\inf})$ but this is not true in general although it is true in many situations, as we will see in the examples below. 
%That is why we formulate the Corollary above in terms of $\lambda_0^-$ and $\lambda_0^+$.  

\rm The results above express the importance of the family of sets $\uK_{\tau_0}, \oK_{\tau_0}$ and the numbers $\lambda_0^{-}$ and $\lambda_0^{+}$ in order to obtain boundedness or unboundedness of solutions in a rather general situation.  Observe that the proposition and the corollary do not say anything if $\lambda \in (\lambda_0^{-},\lambda_0^{+})$.

\begin{figure}[H]
  \centering
    \includegraphics[width=10cm, height=5cm]{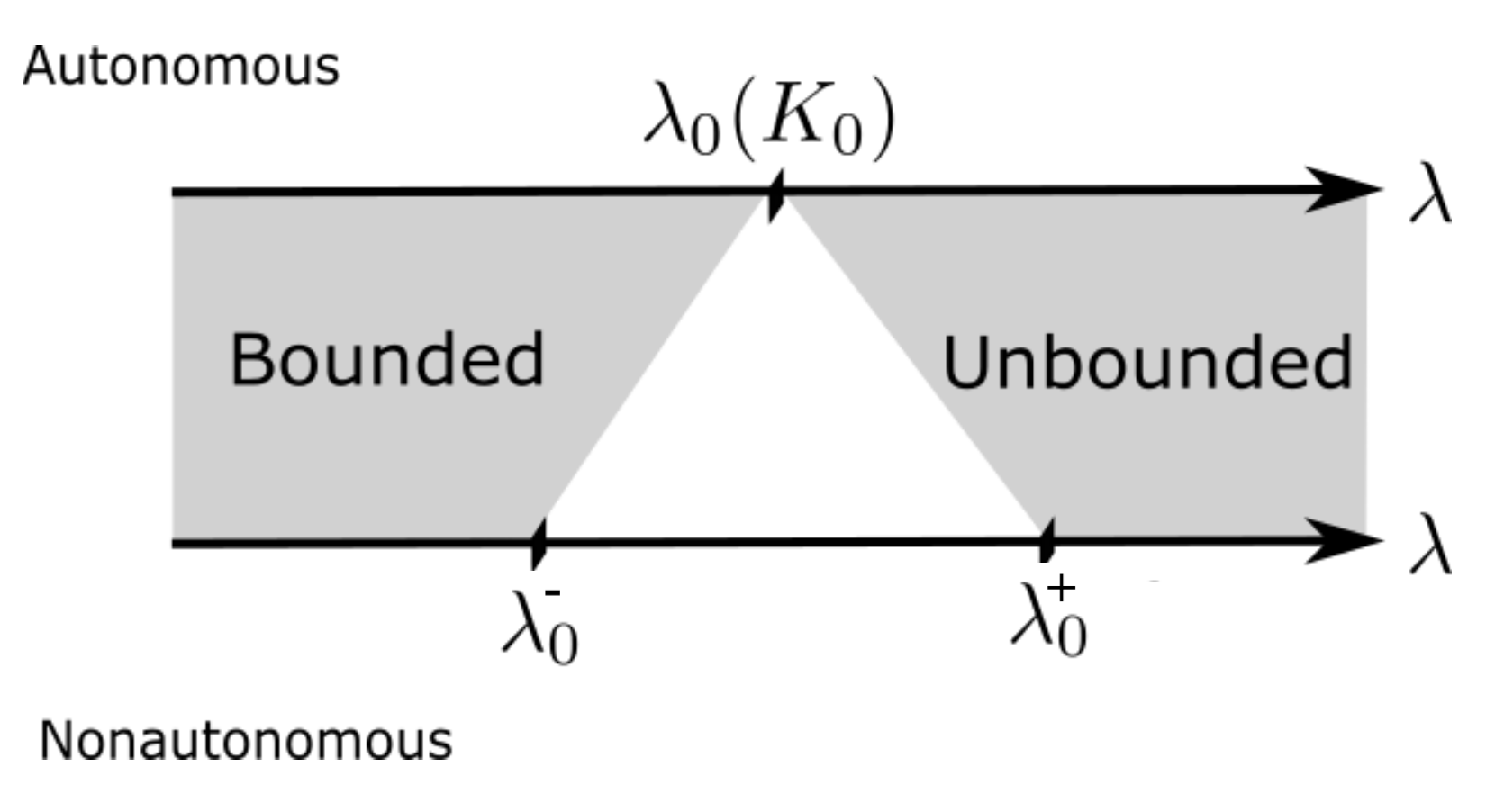}
    \caption{Impact of the time-dependency of $K(t)$ on the variety of scenarios.}\label{autnoaut}
\end{figure}

\rm
\pa Let us see some examples of sets $\uK_{\tau_0}$ and $\oK_{\tau_0}$.

\ex[A Ball of Changing Radius]\label{shrink}  Let $x_0\in\Omega$ and $r_0>0$ such that $\overline{B(x_0,r_0)}\subset \Omega$ and let

\begin{equation*}
n(t,x)=1-\chi_{K(t)}(x), \quad \mbox{for } t > 0,\; x\in \overline{\Omega},
\end{equation*}

\pa where 
\begin{equation*}
K(t)=\overline{B\left(x_0,r(t)\right)}, \; t>0,
\end{equation*}
where $r(t)$ is specified according to different cases.

\medskip

\begin{itemize}

\item[\underline{Case 1}.] Consider $r(t)=r_0(1-\frac{1}{t+1})$. 
Then, $K^{\sup}_{\tau_0}=\overline{B(x_0,r(\tau_0))}$, 
$K^{\inf}_{\tau_0}= \overline{B(x_0,r_0)}$. Moreover,   $\lambda_0(K^{\sup}_{\tau_0})=\lambda_0(\overline{B(x_0,r_0)})(1-\frac{1}{\tau_0+1})^2\to \lambda_0(\overline{B(x_0,r_0)})$ and hence, $\lambda_0^{-}=\lambda_0(\overline{B(x_0,r_0)})$.  Also, since  $K^{\inf}_{\tau_0}$ is independent of $\tau_0$, $\lambda_0^{+}=\lambda_0(K_{\tau_0}^{\inf})$ and if $\lambda < \lambda_0(\overline{B(x_0,r)})$ we have boundedness of solutions and if $\lambda > \lambda_0 (\overline{B(x_0,r)})$ we have unboundedness of solutions.

\item[\underline{Case 2}.]  If $r(t)=\frac{r_0}{t+1}$, that  $K(t)$ is a continuously shrinking ball, then, $\uK_{\tau_0}=B(x_0,\frac{r_0}{\tau_0+1})$ and $\lambda_0^{-}=\lim_{\tau_0\to+\infty} \lambda_0(B(x_0,r_0))(\tau_0+1)^2=+\infty$,  $\oK_{\tau_0}=\{x_0\}$ and $\lambda_0^{+}=+\infty$. Hence solutions are bounded for each $\lambda \in \R$.  

\item[\underline{Case 3}.] Finally, if $r(t)=r_0(1+|\sin(\omega t)|)$ for some $0\neq\omega\in\R$  then, $K^{\sup}_{\tau_0}=\overline{B(x_0,r_0)}$ and $K^{\inf}_{\tau_0}=\overline{B(x_0,2r_0)}$.  Hence if $\lambda < \lambda_0(\overline{B(x_0,2r_0)})$ we have boundedness of solutions and if $\lambda > \lambda_0(\overline{B(x_0,r_0)})$ we have unboundedness of solutions. 
\end{itemize}

\par\medskip 

Observe that in this last case Proposition \ref{prop-boundedness-general} and Corollary \ref{cor-boundedness-general} do not provide any information about the behavior of solutions when $\lambda\in (\lambda_0(\overline{B(x_0,2r_0)}), \lambda_0(\overline{B(x_0,r_0)}))$.

\medskip

\rm We can conclude that the nonautonomous problem presents a broader spectrum of cases than its autonomous version. Also $\lambda$ and $\lambda_0(K(t))$ are no longer  the only elements involved in determining asymptotic behavior of the solutions as, for example, $\lambda_0^{-}$ and $\lambda_0^{+}$ arise naturally in the non--autonomous setting.

\medskip

Let us study other configurations of $K(t)$.

\begin{figure}[H]
  \centering
    \includegraphics[width=6cm]{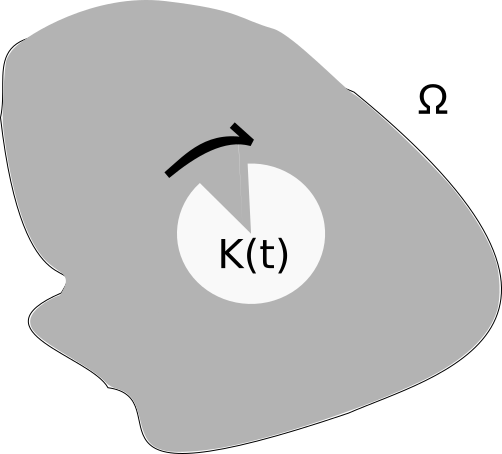}
    \caption{Domain $\Omega$ with a rotating sector.}\label{rotating}
\end{figure}

\ex[Rotating sector]\label{ExKRot} Consider an open domain $\Omega\subset \R^2$ with $\overline{ B(0,r_0)}\subset \Omega$ and let the following set expressed in polar coordinates $$\tilde K=\{ (r,\theta):  0\leq r\leq r_0,\, \theta_0\leq \theta\leq \theta_1\}.$$ 
Define $K(t)$ as the set $\tilde K$, rotating clockwise at certain angular speed $\omega>0$, see Figure \ref{rotating}, that is
$$K(t)=\{(r,\theta):  0\leq r\leq r_0, \,  \theta_0-\omega t\leq \theta\leq \theta_1-\omega t\}.$$ 

In this particular case, $\uK_{\tau_0}$ is the entire circle $\overline{B(0,r_0)}$ as every point $x$ in the circle verifies that $n(t,x)$ vanishes for some $t>\tau_0$ for any $\tau_0>0$. 

%Imagen K rotando

\pa Furthermore, for any $\tau_0>0$, $\oK_{\tau_0}=\{ 0\}$, the center of the circle, as it is the only point where $n$ vanishes for every $t>t_0$. Thus, in this example we have that $\lambda_0(K(t))=\lambda_0(\tilde K)$ is constant and

\begin{equation*}\label{eqRem1cap2}
\lambda_0(\overline{B(0,r_0)})=\lambda_0^{-}<\lambda_0(K(t))<\lambda_0^{+}=\infty , \; \mbox{ for any } t>0.
\end{equation*}

\pa Now, although from the previous proposition we can show that the solutions are bounded for $\lambda < \lambda_0(\overline{B(0,r_0)})$, we cannot decide the behavior for $\lambda >\lambda_0(\overline{B(0,r_0)})$ since we cannot apply Proposition \ref{prop-boundedness-general} to arrive to a conclusion. In fact, as seen later on, whether the solution is bounded or not will also depend on the speed at which $K(t)$ rotates.\rm

\bigskip

\ex[Jumping Sets]\label{Kjump} Let us consider problem (\ref{PreEq1}) and $n(t,x)$ such that $K(t)$ ``jumps'' from a compact set $K_0$ to another compact set $K_1$ periodically:

\begin{equation*}
K(t)=\left\{ \begin{array}{lr}
K_0, & t\in (n t_0, n t_0+T_1], n\in\N,\\
K_1, & t\in (n t_0+T_1, (n+1)t_0], n\in\N,\\
\end{array}\right.
\end{equation*}

\begin{figure}[H]
  \centering
    \includegraphics[width=6cm]{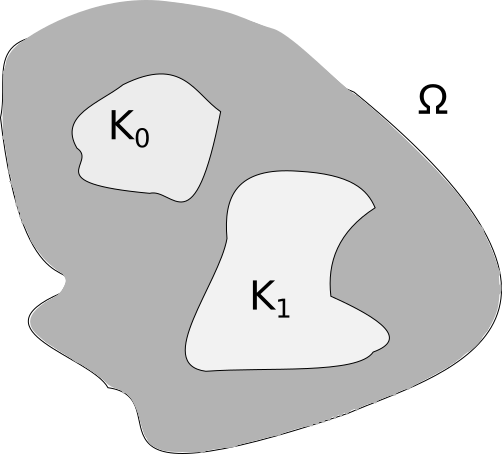}
    \caption{Domain $\Omega$ with ``jumping'' subdomains $K_0$ and $K_1$.}\label{alterK0K1}
\end{figure}

\pa where $0<T_1<t_0$ and $K_0,K_1$ are compact disjoint sets as in Figure \ref{alterK0K1}. This is an example of  special significance as, following the autonomous logic, it is reasonable to expect the existence of unbounded solutions for $\lambda>\lambda_0(K(t))$ for all $t>t_0$. However,  as shown later on, for every $\lambda\in\R$ the solutions remain bounded. The relationship between a jumping $K(t)$ and boundedness will be explored in the following Section \ref{secBound}. \rm

%Imagen Omega y K0 y K1 disconexos

\pa From the previous two examples, it is clear that  the relation between $\lambda$ and $\lambda_0(K(t))$ although important, does not  completely assess the asymptotic behavior of solutions of (\ref{PreEq1}): the geometry of $K(t)$ and how $K(t)$ moves also have to be  factored  in. In fact, we will present several examples of boundedness and unboundedness linked to the geometry of $K(t)$: not only its size, but also taking into consideration how $K(t)$ moves.

\medskip

\section{Boundedness}\label{secBound}
The main objective of this section is to obtain boundedness results that go beyond Proposition \ref{prop-boundedness-general} and Corollary \ref{cor-boundedness-general} in the previous section. In these results, we have established boundedness for $\lambda < \lambda_0^{-}$ and unboundedness for $\lambda >\lambda_0^{+}$.  This section is devoted to explore under which conditions we can obtain boundedness for $\lambda \in (\lambda_0^{-},\lambda_0^{+})$.

Let us start with the particular but important situation where $K(t)=\emptyset$ for some sequence of intervals of times.  We can show:

%%%%%%%%%%%%%%%%%%%%%%%%
\prop\label{PropInter} Let  us assume there exists a sequence of times $t_0<t_1<t_2\ldots<t_i\to +\infty$ and constants  $\eta$, $\Xi >0$  such that \begin{equation}\label{seqT}\left\{\begin{array}{lr}
\left |t_{2i-1}-t_{2i-2}\right|\leq \Xi, & i=1,2,\ldots,\\
\left |t_{2i}-t_{2i-1}\right| > \eta, & i=1,2,\ldots.
\end{array}\right.\end{equation}
Assume also that there exists a $\nu_0>0$ with $n(t,x)\geq \nu_0$ for $t\in (t_{2i-1},t_{2i})$ for all $i\in\N$.  Then, for all $\lambda\in\R$, all solutions of (\ref{PreEq1}) are bounded. Moreover, there exist $\tau>0$ and $\Lambda>0$ independent of $u_0$ such that

\begin{equation*}
\left| u(t,x; t_0,u_0)\right| \leq \Lambda, \quad \forall t\geq\tau, \forall x\in\Omega.
\end{equation*}

\proof Since the boundedness character of solutions is independent of initial conditions due to Proposition \ref{IndepDatIni}  we may take a particular nontrivial initial conditon  $u_0\in L^\infty(\Omega)$.  
\par\medskip First, we  estimate the growth of $u$ in the interval $[t_{2i-2},t_{2i-1})$. In this interval, the only information we have of the function $n$ is that $n(t,x)\geq 0$ but it may vanish in part or all of the domain. We can always bound the growth of $u$ with the linear equation (that is, with $n\equiv 0$). With the notation of Corollary \ref{ExistenciaGlobal}, we have  
\begin{equation*}
\|u(t,\cdot; t_{2i-1}, u(t_{2i-2},\cdot;t_0,u_0))\|_\infty\leq Me^{(\lambda -\lambda_1^\Omega)(t-t_{2i-2})}\|u(t_{2i-2},\cdot; t_0,u_0)\|_\infty, \; t\in [t_{2i-2},t_{2i-1}), i\in\N.
\end{equation*}

\pa Moreover, since $|t_{2i-1}-t_{2i-2}|\leq \Xi$, we have

\begin{equation*}
\|u(t,\cdot; u(t_{2i-2},\cdot;u_0))\|_\infty\leq Me^{\left|\lambda -\lambda_1^\Omega\right|\Xi}\|u(t_{2i-2},\cdot; u_0)\|_\infty, \; t\in [t_{2i-2},t_{2i-1}), i\in\N.
\end{equation*}

\pa Due to the fact that in $[t_{2i-1},t_{2i})$ the nonlinear term is applied to all the domain, we can obtain a supersolution for the time intervals $[t_{2i-1},t_{2i})$ for $i\in\N$ by using the following ODE problem:

\begin{equation}\label{eqODE}
\left\{ 
\begin{array}{lr}
W(t)_t=\lambda W(t) - \nu_0 W(t)^\rho, &  t\in (0,t_{2i}-t_{2i-1}),\\
W(0)=\|u(t_{2i-1},\cdot)\|_{L^\infty(\Omega)},& t=0,
\end{array}
 \right.
\end{equation}

\pa that is, $W(t-t_{2i-1})\geq u(t,x)$ for $t\in (t_{2i-1},t_{2i}]$ for every $i\in\N$ and a.e.  $x\in\Omega$.  As a matter of fact the function $W(t)$ is a solution of

\begin{equation*}\left\{\begin{array}{lr}
W_t(t)-\Delta W(t) =\lambda W(t) - \nu_0 W(t)^\rho\geq \lambda W(t) - n(t,x) W(t)^\rho, &  t\in (0,t_{2i}-t_{2i-1}),\; x\in \Omega,\\
W(t,x)\geq 0,&  t\in (0,t_{2i}-t_{2i-1}),\; x\in \partial\Omega,\\
W(0)=\|u(t_{2i-1},\cdot)\|_{L^\infty(\Omega)},& t=0,\; x\in \Omega,
\end{array}\right.\end{equation*}

\pa which shows that $W$ is a supersolution of $u$. By direct inspection the solution $W$ of (\ref{eqODE}) is:

\begin{equation*}
W(t)=\left[\frac{\nu_0}{\lambda}\left(1-e^{-\lambda (\rho -1)t} \right)+ e^{-\lambda (\rho -1)t} W(0)^{1-\rho}\right]^{-\frac{1}{\rho-1}},
\end{equation*}

\pa which is bounded above by:

\begin{equation}\label{Winf}
W_\infty(t)=\left[\frac{\nu_0}{\lambda}\left(1-e^{-\lambda (\rho -1)t} \right)\right]^{-\frac{1}{\rho-1}}.
\end{equation}

\begin{figure}[H]
  \centering
    \includegraphics[width=10cm]{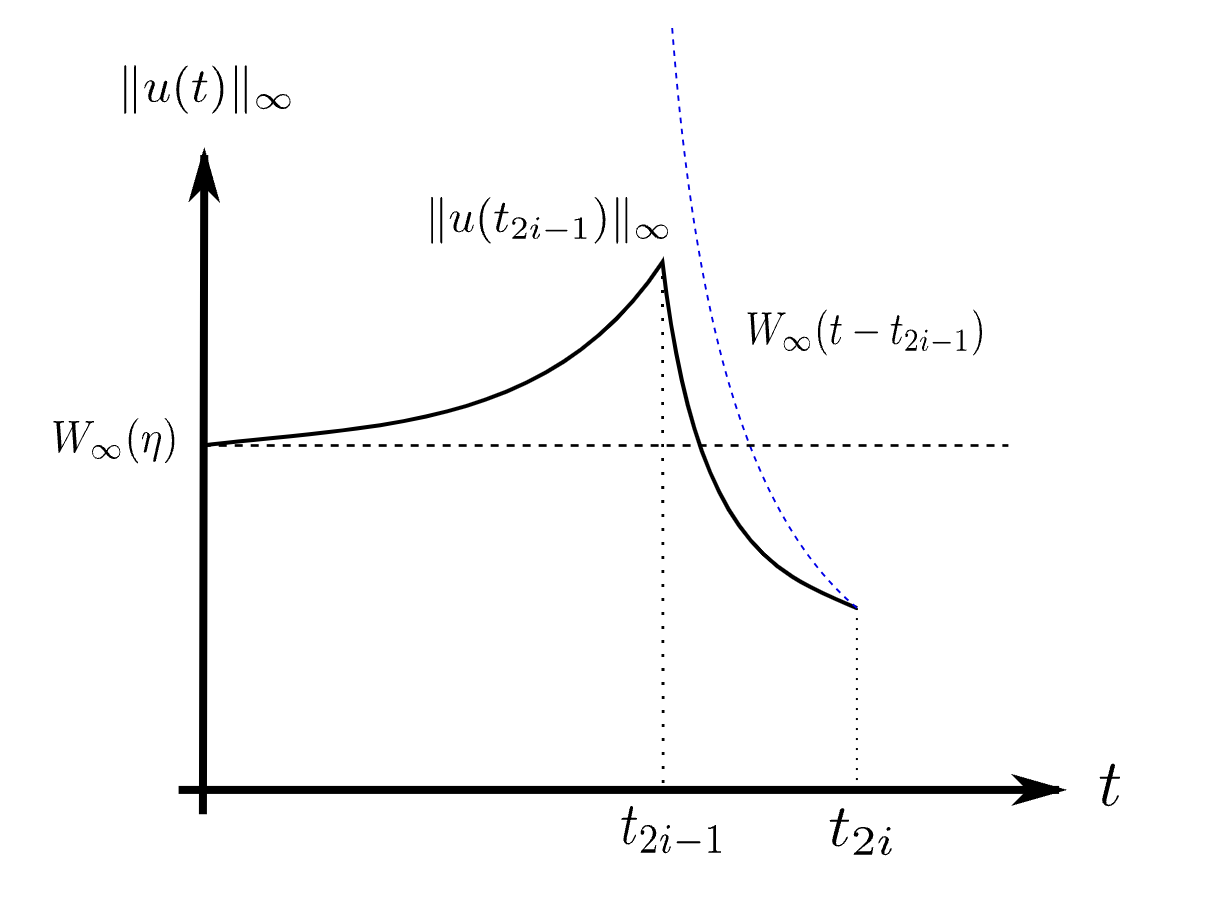}
    \caption{$W_\infty$ bounds $u$ for $t\in (t_{2i-1},t_{2i}]$. }\label{intermitent2}
\end{figure}

The function $W_\infty(t - t_{2i-1})$ is specially important not only because it bounds $u(t,x)$ for $t\in[t_{2i-1},t_{2i})$ for every $i\in\N$ and every $x\in\Omega$ but  it also reaches, for a fixed amount of time $(t_{2i}-t_{2i-1})$, the same threshold $W_\infty(t_{2i}-t_{2i-1})$ wherever the solution $u$ may start at time $t_{2i-1}$. Now, since $|t_{2i}-t_{2i-1}|> \eta$ for every $i\in\N$, we have that $W_\infty(\eta)$ bounds $u(t_{2i})$ for $i\in\N$ and a.e. $x\in\Omega$, see Figure \ref{intermitent2}. In particular,

$$\left| u(t_2,x)\right| \leq W_\infty(\eta),\quad \mbox{a.e.}\; x\in\Omega.$$

%Imagen Omega y K0 y K1

\begin{figure}[H]
  \centering
      \includegraphics[width=10cm]{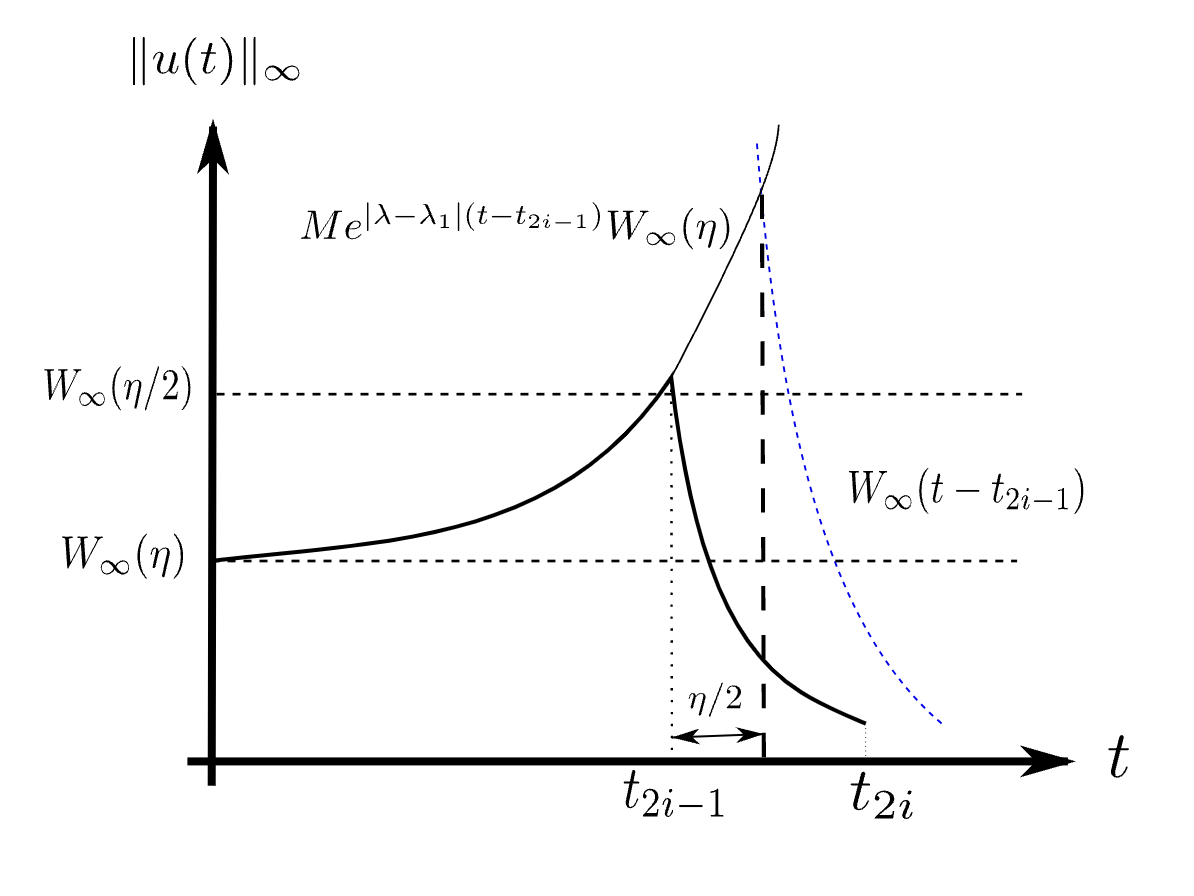}
    \caption{Construction of the upper bound for $u$. }\label{intermitent3}
\end{figure}

\pa Thus,  if we define the function (see Figure \ref{intermitent3}), 

\begin{equation*}
\overline{u}(t):=\left\{ \begin{array}{lr}
M e^{\left|\lambda-\lambda_1^\Omega\right|(t-t_{2i-2})}W_\infty(\eta), & t\in(t_{2i-2},t_{2i-1}+\eta/2],\; i=2,3,\ldots\\
W_\infty(\eta/2), & t\in(t_{2i-1}+\eta/2,t_{2i}],\; i=2,3,\ldots.
\end{array}\right.
\end{equation*}

%Imagen Omega y K0 y K1

 \pa then, we have
 
 \begin{equation*}
\left| u(t,x;u_0) \right| \leq \bar u(t), \quad t\geq t_2,\; x\in\Omega
 \end{equation*}

\pa and in particular

\begin{equation*}
\left| u(t,x;u_0) \right| \leq \Lambda:= \max\limits_{s\in  [t_2,t_4]} \bar u(s), \quad \forall t\geq  t_2.
 \end{equation*}
\pa   Choosing $\tau=t_2$, we conclude the proof of the proposition. \eproof

\re Observe that Proposition \ref{PropInter} applies to all $\lambda \in \R$. In this case, $\lambda_0^+ = \infty$ since $\oK_{\tau_0}=\emptyset$. \rm

\medskip

\subsection{Estimates of the solution locally in space and time}\label{subsecacotlocal}

The effect of having the logistic term acting over the whole domain $\Omega$ during certain intervals of time, allows us to obtain global bounds in space with the use of the function $W_\infty(t)$ (see Proposition \ref{PropInter}). In this subsection, which is inspired by the results from \cite{Du2},   we will obtain bounds of the solutions which are local in space and time under certain conditions on the vanishing logistic term.

Let us start with a partial result that gives us estimates of the solution locally in space but globally in time. In order to accomplish this, we will need to consider the following singular elliptic Dirichlet problem

\begin{equation}\label{z}\left\{\begin{array}{lr}
-\Delta z = \lambda z - \beta z^\rho, & \mbox{ in}\;B(0,a),\\
z=\infty, & \mbox{ on}\;\partial B(0,a).
\end{array}\right.\end{equation}

\pa Problem (\ref{z}) has been studied in numerous articles, see \cite{BandleDiazDiaz,Du0,Du2}. It is known that it has a unique positive radial solution denoted by $z_a(x)$ and this solution is used to construct a supersolution in any ball where the logistic term $n(t,x)$ is strictly positive. In particular, if $\overline{B(x_0,a)} \subset \Omega$ and $n(t,x)\geq \beta > 0$ for all $x\in B(x_0,a)$, $t\geq t_0$ then if $u_0(x) \in L^\infty(\Omega)$ is such that

\begin{equation*}
0 \leq u_0(x) \leq z_a(0)=\inf\limits_{y\in B(0,a)}z_a(y)\leq z_a(y),\quad \forall x\in B(x_0,a),\; \forall y\in B(0,a).
\end{equation*}

\pa Hence,  $z_a(x-x_0)$ is a supersolution of $u(t,x)$, that is,

\begin{equation*}
0 \leq u(t,x) \leq z_a(x-x_0),\quad x\in B(x_0,a),\; t\geq t_0,
\end{equation*}

\pa and therefore, if $C=\sup\limits_{|y|\leq b} z_a(y)$ for $0<b<a$ then 

\begin{equation*}
0 \leq u(t,x) \leq C,\quad x\in B(x_0,b),\; t\geq t_0.
\end{equation*}

\pa Observe that  $C=C(\lambda,\beta,\rho,a,b)$.

\medskip

\pa Now we can show the following result:

\lem\label{acotuK} Let us consider problem (\ref{PreEq1}) with  $n(t,x)$ satisfying {\bf(N)} and let $K\subset \overline{\Omega}$ be a compact set with $K(t)\subset K$ for every $t\geq t_0$ (for instance $K=\uK_{t_0}$).

For each $\delta>0$ small enough there exists $\Lambda=\Lambda(\Omega,\lambda,\rho,n,u_0,\delta)>0$ satisfying

\begin{equation}\label{local1prop}
0\leq u(t,x)\leq \Lambda,\; \mbox{ for } t>t_0, x\in \Omega\setminus\Omega_\delta,
\end{equation}

\pa where $\Omega_\delta =\{ x\in \Omega: d(x,K)<\delta\}$ is an open neighbourhood of $K$.

\proof  Let us consider $\Omega_{\delta/2}$ and define

\begin{equation*}n_{\delta/2}(x):=\left\{\begin{array}{ll}
0,& x\in \overline{\Omega_{\delta/2}},\\
\nu(\delta/2),& x\in \Omega\setminus \overline{\Omega_{\delta/2}},
\end{array}\right.\end{equation*}

\pa where $\nu(\cdot)$ is the function in {\bf{(N)}}. It is clear that $n(t,x)\geq n_{\delta/2}(x),$ $\forall x\in \overline{\Omega_{\delta/2}}$, $t\geq t_0$. Moreover, from {\bf(N)} we have $n(t,x)\geq \nu(d(x,K(t)))$ for every $x\in\Omega\setminus\overline{\Omega_{\delta/2}}$ But since $d(x,K(t))\geq \delta/2$ for $x\in \Omega\setminus\overline{\Omega_{\delta/2}}$ and $\nu$ is increasing, then we have $n(t,x)\geq \nu(\delta/2)>0$. Hence, $n(t,x)\geq n_{\delta/2}(x)$ $\forall x\in \overline{\Omega}$ and $\forall t\geq t_0$ and by comparison $ u(t,x;u_0)\leq \hat u(t,x,;u_0)$ where $\hat u$ is the solution of \eqref{eq2intro} with $n(x)=n_{\delta/2}(x)$. 

But estimates (\ref{local1prop}) were shown in \cite{AnibalArrietaPardo} for $\hat u$ (since $\hat u$ solves an autonomous problem) as it has been shown above.\eproof

\medskip

\begin{figure}[H]
  \centering
    \includegraphics[width=8cm, height=5cm]{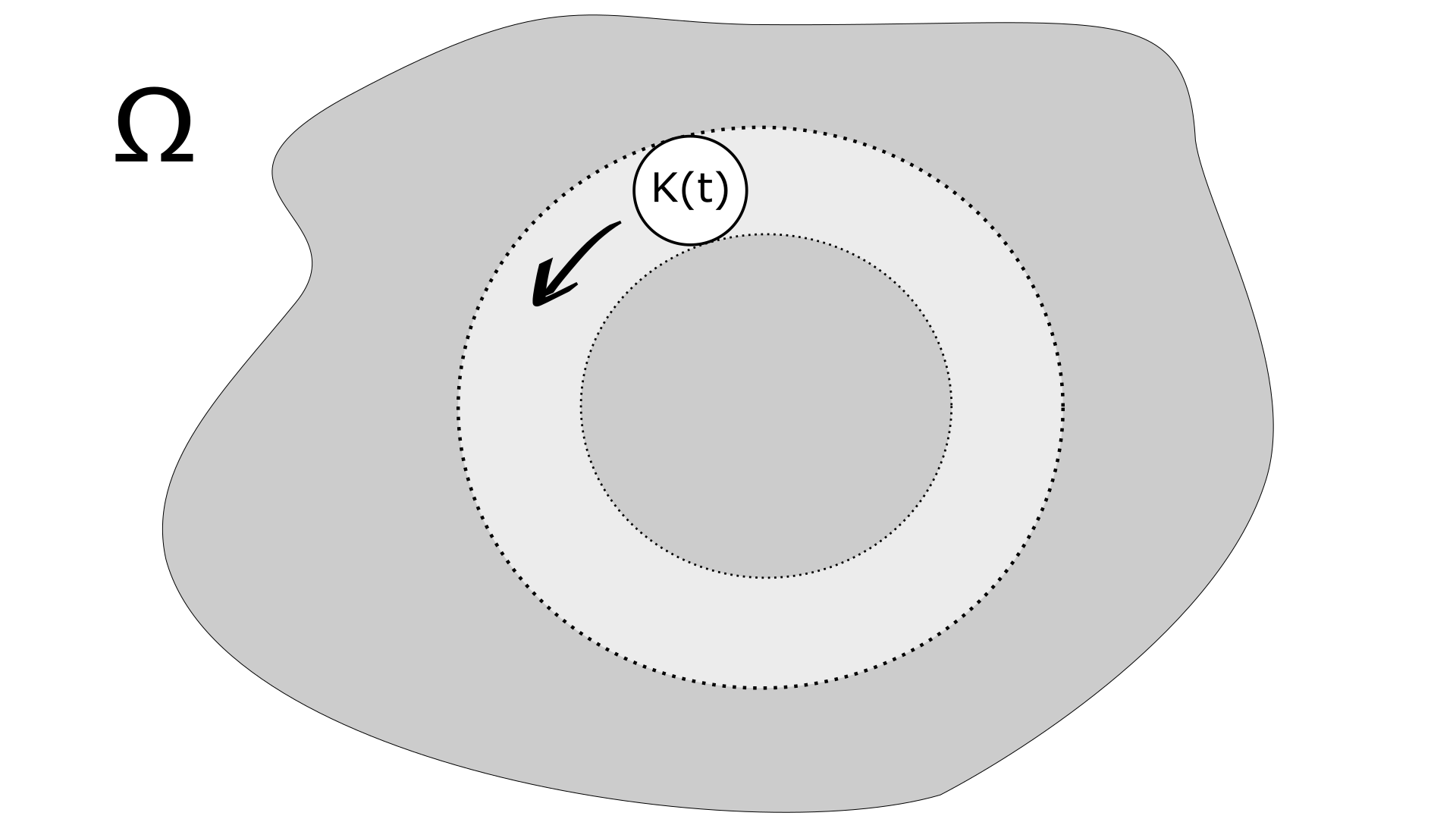}
    \caption{A continuously moving $K(t)$ along the ring $\bigcup\limits_{t\geq t_0}K(s)$.}
    \label{K(t)-moving ring}
\end{figure}

\pa Note that the bounds obtained in the previous Lemma \ref{acotuK} apply only to points outside $\uK$. But in many situations we may have points $x_0\in\Omega$ with $\overline{B(x_0,a)}\subset \Omega\setminus K(t)$ for $t\in (t_1,t_2)$ while $x_0\in K(t)$ for other times. For instance, consider  $K(t)$ is a continuously moving ball within a ring for which $\uK$ is the ring, see Figure \ref{K(t)-moving ring}. If $x\in \uK$ we have times when $x\in K(t)$, but we have other times where $x\in \Omega\setminus K(t)$. It is very relevant to study the behavior of the solution in these points.

%Imagen K rotando

\pa We have the following:

\prop\label{acotlocal1} Let  $u$ be the solution of (\ref{PreEq1}) for $0\leq u_0\in L^\infty(\Omega)$. Let $x_0\in\Omega$ and $a>0$ be such that $B(x_0,a)\subset \Omega$ and let $n$ satisfy 

\begin{equation*} 
 n(t,x)\geq \beta>0,\; t\in(t_1,t_2),\; x\in B(x_0,a),
\end{equation*}

\pa for some $\beta>0$ and $t_1<t_2$. Then, for any $\tau_0\in (0,t_2-t_1)$ and any $r\in (0,a)$ there exists $\Lambda=\Lambda(\lambda,\beta,\rho,\tau_0,r,a)$ such that 

\begin{equation*}
0\leq u(t,x) \leq \Lambda,\quad\mbox{for } t\in [t_1+\tau_0,t_2],\; x\in B(x_0,r).
\end{equation*}

\proof Let us consider $U^*(t,x)=z_{a,x_0}(x)+W_\infty(t-t_1)$ for $t\in (t_1,t_2)$ and $x\in B(x_0,a)$ where $z_{a,x_0}(x)=z_a(x-x_0)$ and $z_a$ is the unique positive solution of (\ref{z}) and $W_\infty$ is defined in (\ref{Winf}). Then, since $z_a,W_\infty >0$ and $\rho>1$, we have
\begin{equation*}\begin{array}{rl}
U^*_t-\Delta U^* & = {W_\infty}_t-\Delta z_{a,x_0}=\lambda W_\infty - \beta W_\infty^\rho +\lambda z_{a,x_0} - \beta z_{a,x_0}^\rho\\
\\
\; & = \lambda U^* - \beta\left( W_\infty^\rho+ z_{a,x_0}^\rho\right) = \lambda U^* -\beta(W_\infty+ z_{a,x_0})^\rho\left( \left(\frac{ W_\infty}{W_\infty+ z_{a,x_0}}\right)^\rho+\left(\frac{ z_{a,x_0}}{W_\infty+ z_{a,x_0}}\right)^\rho  \right)\\
\\
\; & \geq \lambda U^* -\beta(W_\infty+ z_{a,x_0})^\rho\left( \left(\frac{ W_\infty}{W_\infty+ z_{a,x_0}}\right)+\left(\frac{ z_{a,x_0}}{W_\infty+ z_{a,x_0}}\right)  \right)= \lambda U^* -\beta(U^*)^\rho.
\end{array}\end{equation*}

\pa Moreover $U^*(t,x)=\infty$ for $t\in (t_1,t_2)$, $x\in\partial B(x_0,a)$ and for $t=t_1$, $x\in B(x_0,a)$. Therefore, $U^*$ is a supersolution for $u$ in $(t_1,t_2)\times B(x_0,a)$ for any initial condition $u_0\in L^\infty(\Omega)$. Notice that for any $\tau_0\in(0,t_2-t_1)$ and any $r\in (0,a)$ if $ t\in[t_1+\tau_0,t_2),\; x\in B\left(x_0,r\right)$ then

\begin{equation*}
0\leq u(t,x) \leq U^*(t,x) \leq \|W_\infty(\tau_0)\|_\infty+\|z_{a,x_0}(\cdot)\|_{L^\infty\left(B\left(0,r\right)\right)}<\infty.
\end{equation*}

\pa Which concludes the proof, choosing $\Lambda=\|W_\infty(\tau_0)\|_\infty+\|z_{a,x_0}(\cdot)\|_{L^\infty\left(B\left(0,r\right)\right)}$.\eproof

\re  Other bounds for $U^*$ when $x\rightarrow \partial B(x_0,a)$ or, $t\rightarrow t_2$ for $\beta=\beta(t)$ continuous, $\beta(t)>0$ for every $t\in (t_1,t_2)$ and $\beta(t_2)=0$  can be seen in \cite{Du1}.\rm

\medskip

With Proposition \ref{acotlocal1} it is not difficult to obtain bounds not only in balls but in other subsets at a positive distance from the sets $K(t)$. 

\cor\label{acotlocSubdom} Let $D$ be a subset of $\Omega$,  $\tau_0>0$, and $\delta>0$. If $d\left(D,\bigcup\limits_{t\in[t_1,t_2]} K(t)\right)>\delta>0$  for any $t_1,t_2\in \R$ with $t_2-t_1>\tau>0$, then, there exists $\Lambda=\Lambda(\lambda,\rho,\tau, \nu,\delta)$ such that 

\begin{equation*}
0\leq u(t,x) \leq \Lambda,\quad \forall t\in [t_1+\tau,t_2],\; \forall x\in D.
\end{equation*}

\proof Let us take the open set $\mathcal{U}$

$$\mathcal{U}:=\bigcup\limits_{z\in D}B\left(z,\frac{\delta}{2}\right)\supset D.$$

\pa Since $D$ is at least a distance $\delta$ away from $K(t)$ for $t\in [t_1,t_2]$ we know that $\mathcal{U}$ is at least a distance $\delta/2$ away from $K(t)$ for the same time interval. Therefore,

$$n(t,x)\geq \nu(d(x,K(t))\geq \nu\left(\frac{\delta}{2}\right)>0,\quad \forall t\in[t_1,t_2], \; x\in \mathcal{U}.$$

\pa Moreover, applying Proposition \ref{acotlocal1} for every ball $B(z,\delta/2)$ with $z\in D$ and using the continuity of $u$ in $t$ we obtain that there exists a constant $\Lambda=\Lambda(\lambda,\rho,\tau, \nu(\delta/2))$ such that 

\begin{equation*}
0\leq u(t,x) \leq \Lambda,\quad \forall t\in [t_1+\tau,t_2],\; \forall x\in B\left(z,\frac{\delta}{4}\right),\; \forall z\in D,
\end{equation*}

\pa which concludes the proof. \eproof

\subsection{Boundedness for $\lambda < \lambda_0(K(t)), \forall t$}\label{comienzoPert}

In this subsection we obtain one of the main results of this section which consists of a boundedness result for the case where, roughly speaking, $\lambda < \lambda_0(K(t))$ $\forall t\geq t_0$. As a matter of fact, we will show the following important result:

\teo\label{TeoAcot1} Assume the function $n(\cdot)$ satisfies {\bf(N)}. Assume also that there exist  $ \tau_0>0$, $\delta>0$ and a family of smooth bounded open sets $\left\{\Omega_\delta(t)\right\}_{t\geq t_0}$ with $\Omega_\delta(t)\subset \Omega$ such that for each $t\geq t_0+\tau_0$

\begin{equation}\label{condAcotLambdaSmall}
\bigcup\limits_{s\in[t-\tau_0,t+\tau_0]} K(s) \subset \Omega_\delta (t),  \qquad  \mbox{ and } \qquad  d\left(\bigcup\limits_{s\in [t-\tau_0,t+\tau_0]} K(s), \,\, \partial \Omega_\delta(t)\right)\geq \delta.
\end{equation}

\pa Then if $\lambda<\lambda_1\equiv \inf\limits_{t\geq t_0}\lambda_1(\Omega_\delta(t))$ %and $\tau_0>\omega\equiv\sup\{ |\Omega_\delta(t))|^{2/N}/4\pi\,: \,\, t\geq t_0\}$, 
we have that all solutions of (\ref{PreEq1}) are bounded.  Moreover, the solutions are uniformly asymptotically bounded in $L^\infty$. 

\re Condition (\ref{condAcotLambdaSmall}) implies that $K(t)$ is at a distance at least $\delta>0$ from $\partial \Omega$. See Corollary \ref{corExtDomAcot1} below to extend this result to the case where $K(t)\cap \partial\Omega \neq \emptyset$.

\proof %[Proof(of Theorem \ref{TeoAcot1})]

\begin{figure}[H] 
  \centering
    \includegraphics[width=10cm]{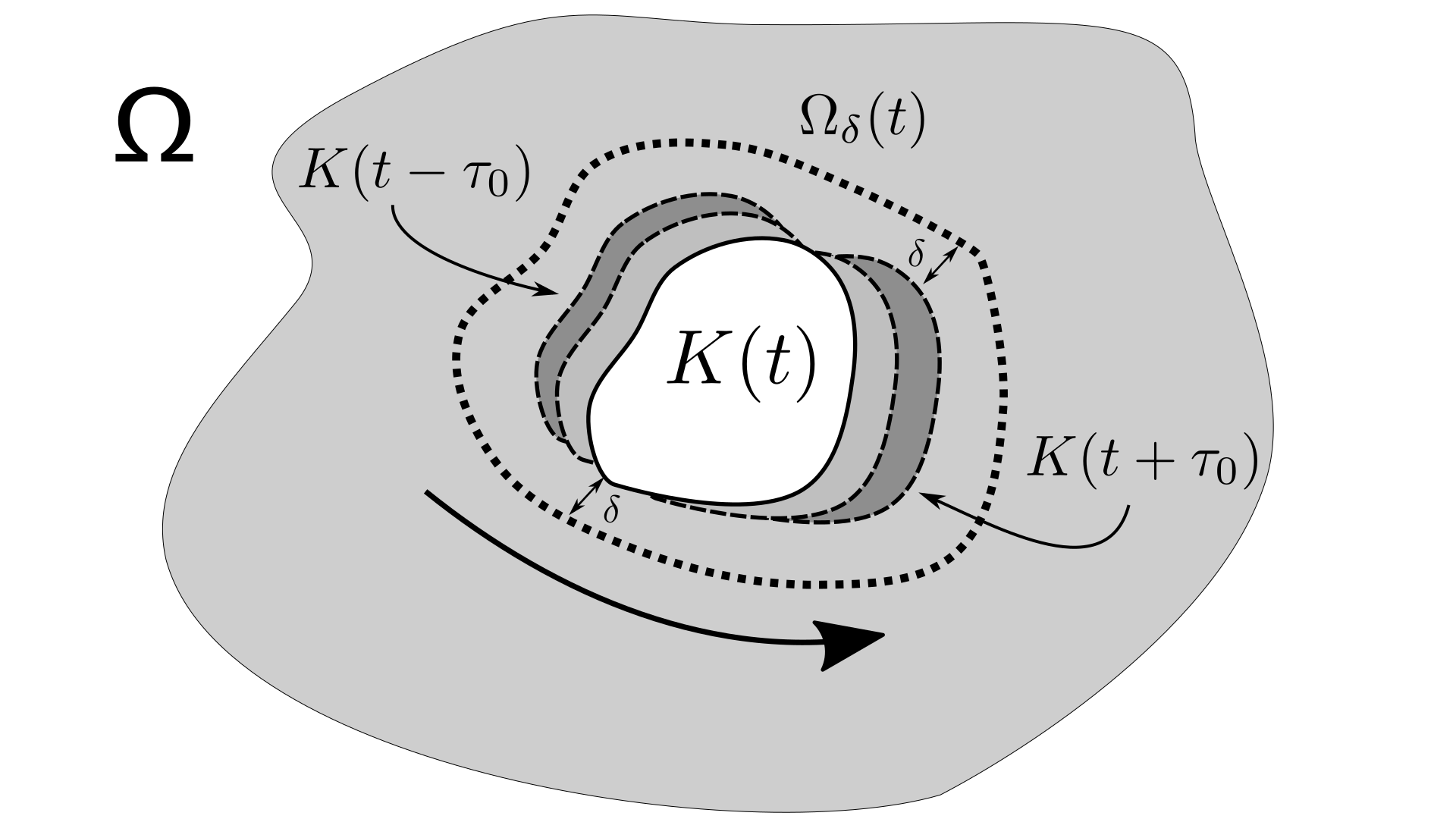}
        \caption{$\Omega_\delta(t)$ contains $\displaystyle\bigcup_{s\in [t-\tau_0,t+\tau_0]} K(s)$ and it lies at a distance $\delta$ from $\partial\Omega_\delta(t)$.}\label{figuraD}
\end{figure}

Let us define $t_1 = t_0+\tau_0$, and in general, $t_{i+1} = t_i +\tau_0$ for $i\geq 1$, $\lambda_1:=\inf\limits_{t\geq t_0} \lambda_1(\Omega_\delta(t))$ and let $D_i$ be a $\delta/2$--neighbourhood of $\Omega\setminus \Omega_\delta(t_i)$, that is, $D_i=\{ x\in \Omega: d(x,\Omega\setminus\Omega_\delta(t_i))< \delta/2\}$ for $i\geq 1$, see Figure \ref{figuraD}.  Now, from (\ref{condAcotLambdaSmall}) we have that for any $x\in D_i$ and $i\geq 1$

 $$d\left(x,\bigcup\limits_{s\in[t_i-\tau_0,t_i+\tau_0]}K(s)\right)\geq \delta/2.$$

\pa Applying Corollary \ref{acotlocSubdom} to $D_i$ for $i\geq 1$ we obtain that there exists a bound $\Lambda=\Lambda(\lambda,\rho,\tau_0,\nu,\delta)$ such that $u(t,x)\leq \Lambda$ for $t\in [t_i,t_i+\tau_0]$ and $x\in D_{i}$.  In particular, since $\Omega \setminus \Omega_\delta(t_i) \subset D_{i}$ and $\partial\Omega_\delta(t_i) \subset D_{i}$,  constant $\Lambda$ is independent of $i$ and satisfies
\begin{equation}\label{LandaNormal}
u(t,x) \leq \Lambda, \qquad \mbox{ for }t\in [t_i,t_{i+1}],\; \mbox{ a.e. } x\in \Omega\setminus \Omega_\delta(t_i),\;  i\geq 1.
\end{equation}

\pa In particular we have $u(t,x) \leq \Lambda$,  for $t\in [t_i,t_{i+1}]$, a.e. $x\in \partial \Omega_\delta(t_i)$ and any $i\geq 1$.

\medskip

 Let us prove that  we have
 
 \begin{equation}\label{first-L2-estimate}
 \|u(t_{i+1},\cdot)\|_{L^2(\Omega_\delta(t_i))}\leq \frac{\lambda|\Omega|^{1/2}}{\lambda-\lambda_1^{\Omega}}\Lambda+e^{(\lambda-\lambda_1)\tau_0}\|u(t_{i},\cdot)\|_{L^2(\Omega_\delta(t_{i}))}, \quad i=1,2,\ldots
 \end{equation}

 \pa  To do so, let us consider the problem

\begin{equation}\label{EqCircling1}\left\{\begin{array}{lr}
v_t-\Delta v = \lambda v , & x\in\Omega_\delta(t_i), \; t\in (t_i,t_{i+1}),\\
v = \Lambda, & x\in\partial\Omega_\delta(t_i), \; t\in (t_i,t_{i+1}),\\
v(t_i,x)= u(t_i,x), & x\in\Omega_\delta(t_i).
\end{array}\right.\end{equation}

By comparison, we have 

$$v(t,x)\geq u(t,x), \qquad t\in [t_i,t_{i+1}], \mbox{ a.e. }x\in\Omega_\delta(t_i).$$
\begin{figure}[H]
  \centering
    \includegraphics[width=10cm]{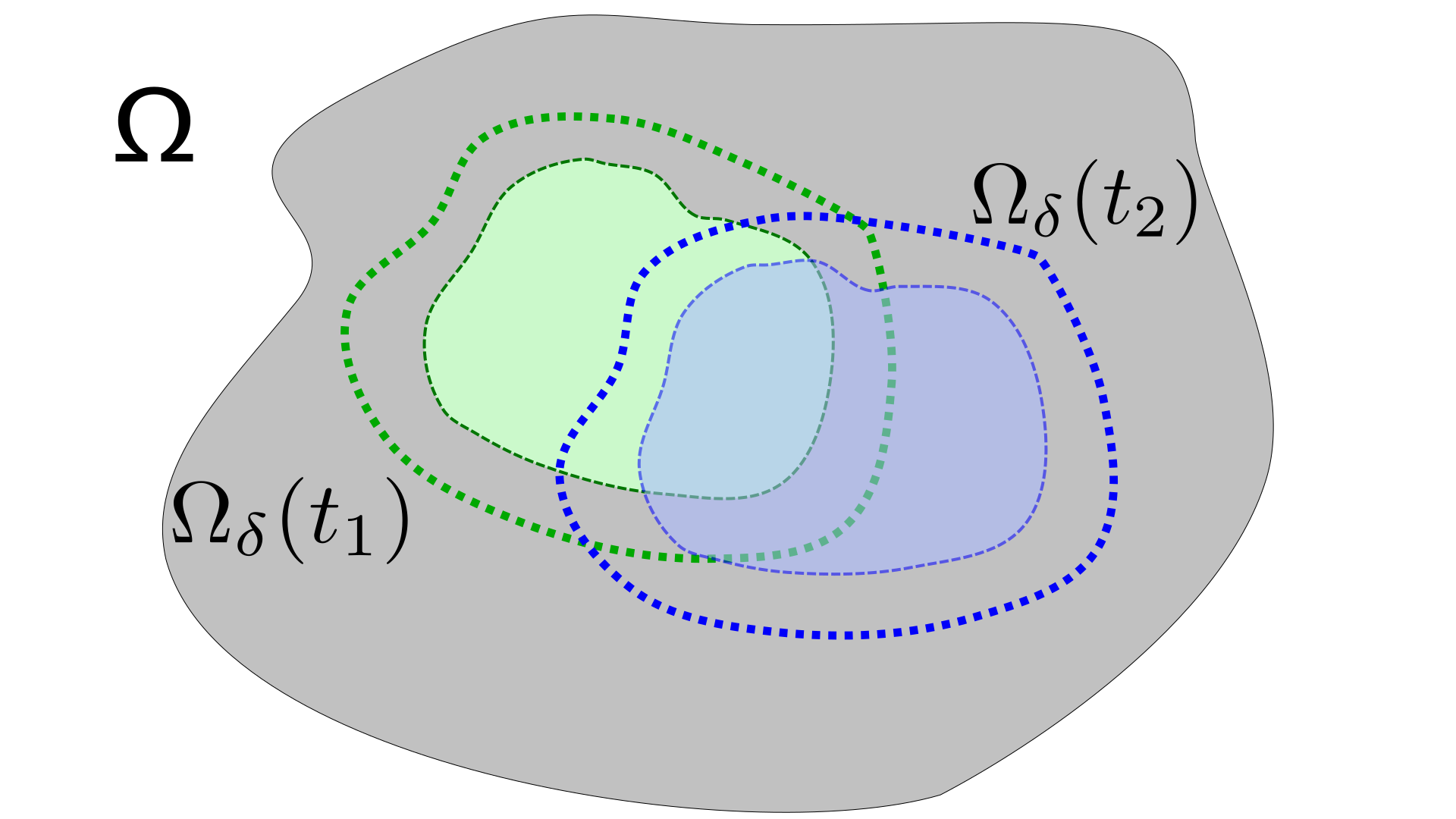}
        \caption{Neighbourhoods $\Omega_\delta(t_1)$ and $\Omega_\delta(t_2)$ (dotted lines). }
\end{figure}

\pa Now, $v$ can be split  as $v=v_1+v_2$ where

\begin{equation}\label{EqV1}\left\{\begin{array}{lr}
-\Delta v_1 = \lambda v_1 , & x\in\Omega_\delta(t_i),\\
v_1= \Lambda, & x\in\partial\Omega_\delta(t_i),
\end{array}\right.\end{equation}

\begin{equation}\label{eqV2teo1pert}
\left\{\begin{array}{lr}
v_{2_t} -\Delta v_2= \lambda v_2 , & x\in\Omega_\delta(t_i), \; t\in (t_i,t_{i+1}),\\
v_2=0, & x\in \partial\Omega_\delta(t_i), \; t\in (t_i,t_{i+1}),\\
v_2(t_i,x)= u(t_i,x)-v_1(x), & x\in\Omega_\delta(t_i).
\end{array}\right.\end{equation}

\pa Observe that since $\lambda<\lambda_1\leq \lambda_1(\Omega_\delta(t_1))$, $v_1$ is well defined. Moreover, standard estimates show that 
$$\|v_1\|_{L^2(\Omega_{\delta}(t_i))}\leq \frac{\lambda |\Omega|^{1/2}}{\lambda-\lambda_1(\Omega)}\Lambda.$$

Since,  by the maximum principle, $v_1(x)\geq 0$ we have that  standard $L^2$-estimates applied to \eqref{eqV2teo1pert} show 
$$\|v_2(t_{i+1},\cdot)\|_{L^2(\Omega_\delta(t_i))}\leq e^{(\lambda-\lambda_1)\tau_0}\|v_2(t_{i},\cdot)\|_{L^2(\Omega_\delta(t_i))}.$$

\pa Hence, since $0\leq u(t_{i+1},x)\leq v_1(x)+v_2(t_{i+1},x)$ from the two previous estimates we obtain \eqref{first-L2-estimate}.

\par\medskip  But observe now that if $x\in \Omega_\delta(t_{i+1})\setminus \Omega_\delta(t_{i})$ we have $u(t_{i+1},x)\leq \Lambda$. Hence
$$\|u(t_{i+1},\cdot)\|_{L^2(\Omega_\delta(t_{i+1}))}^2\leq \|u(t_{i+1},\cdot)\|_{L^2(\Omega_\delta(t_{i+1})\cap \Omega_\delta(t_{i})   )}^2+\Lambda^2|\Omega_\delta(t_{i+1})\setminus \Omega_\delta(t_{i})|\leq \|u(t_{i+1},\cdot)\|_{L^2(\Omega_\delta(t_{i}))}^2+\Lambda^2|\Omega|$$
and therefore, plugging in \eqref{first-L2-estimate} in this last expression we get

\begin{equation}\label{second-L2-estimate}
\|u(t_{i+1},\cdot)\|_{L^2(\Omega_\delta(t_{i+1}))}\leq \frac{2\lambda -\lambda_1}{\lambda-\lambda_1}|\Omega|^{1/2} \Lambda+e^{(\lambda-\lambda_1)\tau_0}\|u(t_{i},\cdot)\|_{L^2(\Omega_\delta(t_{i}))}, \quad i=1,2,\ldots
\end{equation}
If we denote by $z_i=\|u(t_{i},\cdot)\|_{L^2(\Omega_\delta(t_{i}))}$,  $A=\displaystyle \frac{2\lambda -\lambda_1}{\lambda-\lambda_1}|\Omega|^{1/2} \Lambda$ and $\rho=e^{(\lambda-\lambda_1)\tau_0}$, we have $z_{i+1}\leq A+\rho z_i$, $i=1,2,\ldots$ and iterating this inequality we obtain
$$\|u(t_{i+1},\cdot)\|_{L^2(\Omega_\delta(t_{i+1}))}\leq A(1+\rho+\ldots+\rho^{i-1})+\rho^{i}\|u(t_{1},\cdot)\|_{L^2(\Omega_\delta(t_{1}))}\leq \frac{A}{1-\rho}+\rho^i\|u(t_{1},\cdot)\|_{L^2(\Omega_\delta(t_{1}))}.$$

Now since $t_1=t_0+\tau_0$, standard $L^2$-estimates in $\Omega$ show that

$$\|u(t_{1},\cdot)\|_{L^2(\Omega_\delta(t_{1}))}\leq \|u(t_{0}+\tau_0,\cdot)\|_{L^2(\Omega)}\leq e^{(\lambda-\lambda_1^\Omega)\tau_0}\|u_0\|_{L^2(\Omega)},$$
and choosing $i_0\geq 1$ large enough so that $\rho^{i_0}e^{(\lambda-\lambda_1^\Omega)\tau_0}\|u_0\|_{L^2(\Omega)}\leq 1$, we get 
$$\|u(t_{i+1},\cdot)\|_{L^2(\Omega_\delta(t_{i+1}))}\leq  \frac{A}{1-\rho}+1,\qquad i\geq i_0.$$

But, since $u(t_{i+1},x)\leq \Lambda$ for $x\in \Omega\setminus\Omega_\delta(t_{i+1})$ and taking into account the definition of $A$ and $\rho$, we have that if we define
$$R=\Lambda|\Omega|^{1/2}+ \frac{2\lambda -\lambda_1}{\lambda-\lambda_1}|\Omega|^{1/2} \Lambda \frac{1}{1-e^{(\lambda-\lambda_1)\tau_0}}+1,$$
then
$$\|u(t_{i+1},\cdot)\|_{L^2(\Omega)}\leq R,\qquad i\geq i_0.$$
Using now $L^2-L^\infty$-estimates, we get, 
$$\|u(t,\cdot)\|_{L^\infty(\Omega)}\leq CR, \qquad t\in [t_{i+2},t_{i+3}], \quad i\geq i_0.$$
This shows the result. \endproof

\re\label{RecondAcotLambdaSmall} Conditions (\ref{condAcotLambdaSmall}) can be easily met, for example, whenever $K(t)=\gamma(t)+R(t) K_0$, where $\gamma(t)$ is a $\mathcal{C}^1$ curve in $\Omega$ with $\sup_{t\geq t_0}|\gamma'(t)|<\infty$, $R(t)$ is a rigid motion for all $t\geq t_0$ and is continuous with respect to $t\geq t_0$ with the condition that $\gamma(t)+R(t) K_0\subset \Omega$ $\forall t\geq t_0$, and $K_0$ is a compact subset. If $\lambda<\lambda_0(K_0)$ and $K_0$ can be seen as the intersection of a family of nested open neighbourhoods of $K_0$, then there exists an open neighbourhood $\Omega_\delta$ such that $\lambda < \lambda_1(\Omega_\delta)$ and $d(\partial\Omega_\delta, K_0)>2\delta$ for some $\delta>0$.

\pa Since $K(t)$ moves continuously, for any positive distance to the boundary of $\Omega_\delta$, let us consider $\delta$, there exists a time $\tau$ such that 

$$\bigcup\limits_{[t-\tau,t+\tau]} K(s) \subset \Omega_\delta (t),  \qquad  \mbox{ and } \qquad  d\left(\bigcup\limits_{[t-\tau,t+\tau]} K(s), \partial \Omega_\delta(t)\right)\geq \delta.$$

\pa Finally, since there is a maximum speed at which it travels, there is a minimum time $\tau_0$ which verifies all the above.\rm

\bigskip

Finally, as we mentioned above, we can consider the case where $K(t)$ touches $\partial \Omega$.

\cor\label{corExtDomAcot1} Let $R>0$ such that $\overline{\Omega}\subset B(0,R)$. Let us assume $n$ satisfies {\bf{(N)}} and there exist $\tau_0>0,\delta>0$ and a family of smooth bounded open sets $\{\Omega_\delta(t)\}_{t\geq t_0}$ with $\Omega_\delta(t) \subset B(0,R)$ such that for each $t\geq t_0+\tau_0$,  (\ref{condAcotLambdaSmall}) in Theorem \ref{TeoAcot1} is satisfied.

\pa Then, if $\lambda<\inf\limits_{t\geq t_0}\lambda_1(\Omega_\delta(t))$ and $\tau_0>\omega$, we have that all solutions of (\ref{PreEq1}) are uniformly asymptotically bounded in $L^\infty$.

\proof If we define

$$\tilde n(t,x) := \left\{\begin{array}{lr}
n(t,x), & t\geq t_0, \; x\in\overline{\Omega},\\
\nu (2R), & t\geq t_0, \; x\in B(0,R)\setminus\overline{\Omega},
\end{array}\right.$$

\pa we have $\tilde n$ satiesfies condition {\bf{(N)}} also with the same function $\nu$. Moreover, by comparison,  the solution $u(t,t_0, u_0, \tilde n, B(0,R))$ is a supersolution of $u(t,t_0, u_0, n, \Omega)$. But we can apply Theorem \ref{TeoAcot1} to $\tilde u$, obtaining boundedness for $\tilde u$ and therefore for $u$.\eproof

\subsection{Boundedness for a jumping $K(t)$}\label{jumpingsection}

This final boundedness result characterizes boundedness with a discontinuity in $K(t)$ independently of the size of $\lambda$. It is an interesting result since it does not consider the linear growth occurring in the sets $K(t)$ and its relationship with $\lambda$: it only focuses on the fact that the sets  $K(t)$ jump at certain instants of time. 

As a prototype of this situation, see Example \ref{Kjump} and Figure  \ref{ejemploKjump45} below.

%Imagen K rotando
\begin{figure}[H]
  \centering
    \includegraphics[width=9cm]{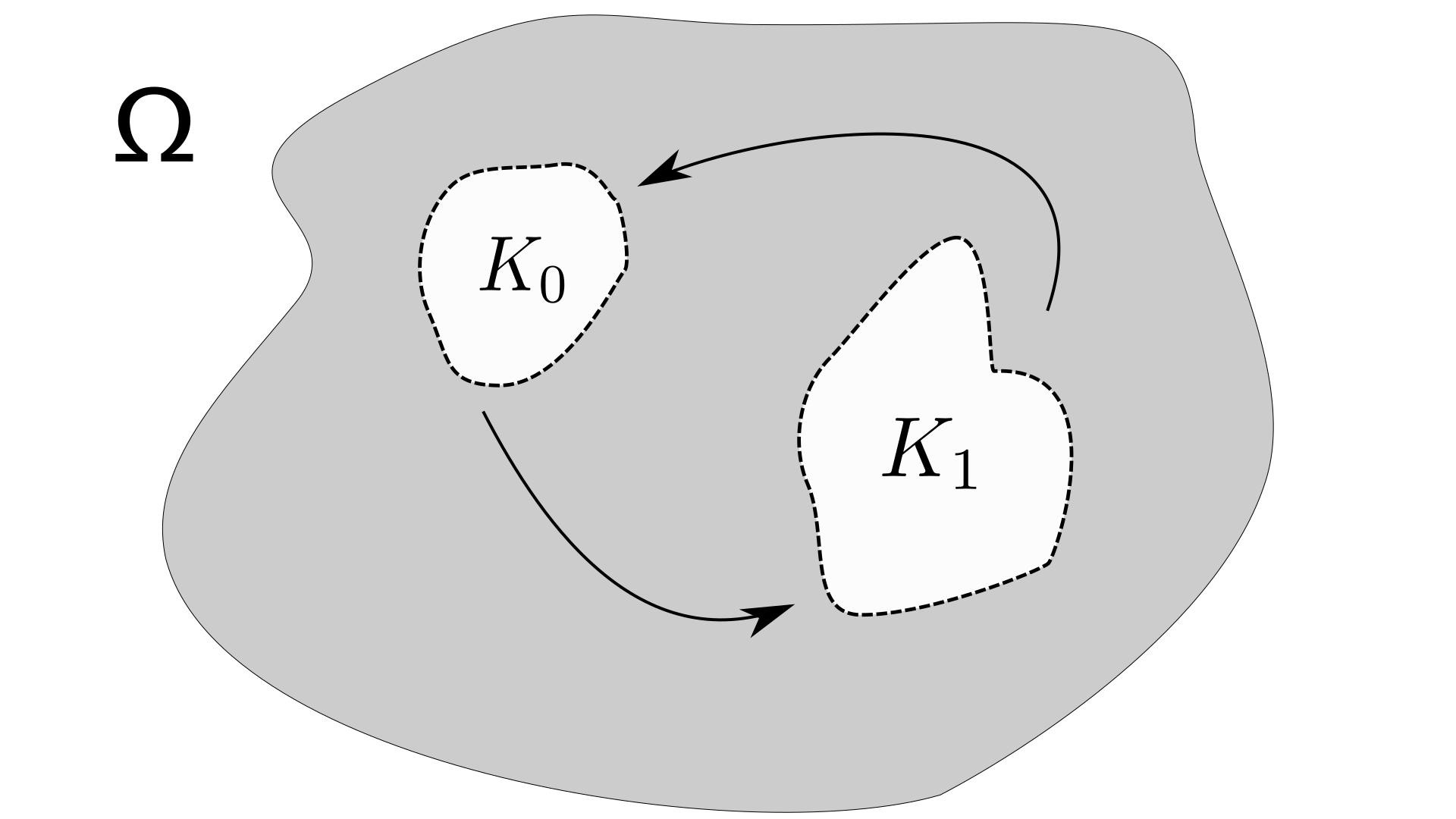}
    \caption{A jumping $K(t)$ between compact subsets $K_0$ and $K_1$.}\label{ejemploKjump45}
\end{figure}

\pa That is, by ``jumps'' or discontinuity in $K(t)$ we refer to a sequence of times $(t_j)_{j\in\N}$ with $t_1<t_2<t_3< \ldots \rightarrow \infty$ where the moving set $K(t)$ jumps at each $t_j$. As a matter of fact we have:

\teo\label{TeoAcot} Let $u$ be the solution of (\ref{PreEq1}), $\lambda\in\R$ and let $n$ satisfying {\bf(N)}. Moreover, let us assume there exists an unbounded increasing sequence $(t_j)_{j\in\N}$ and $\tau_0,\delta,\Xi>0$  such that

\begin{equation}\label{condAcotGen}
\left\{
\begin{array}{l}
2\tau_0\leq \inf\limits_{j\in\N} |t_{j+1}-t_j|\leq \sup\limits_{j\in\N} |t_{j+1}-t_j| =\Xi < \infty \; , \\ \\\qquad d\left(\bigcup\limits_{(t_j-\tau_0,t_j)} K(s),\bigcup\limits_{(t_j,t_j+\tau_0)} K(s)\right) > \delta . \\
\end{array}
\right.
\end{equation}

\pa Then, all solutions of (\ref{PreEq1}) are bounded.

\proof In this proof we will use that $K(t)$ is bounded away from the boundary of $\Omega$. Corollary \ref{corExtDomAcot4} below, extends this result even when $K(t)\cap\partial\Omega\neq \emptyset$.

Let us start showing that  we have a constant $\Lambda^*>0$  such that for any initial condition we can obtain the bounds:  
\begin{equation}\label{claim}
\|u(t_j+\tau_0,\cdot)\|_{L^\infty(\Omega)}\leq \Lambda^*,\quad \hbox{ for }j=1,2\ldots
\end{equation}

Let us consider for each $j\in\N$ the set
$$  \Omega^+_{\delta}(t_j):= \left\{x\in \Omega: d\left(x, \bigcup\limits_{(t_j,t_j+\tau_0)} K(s)\right)<\frac{\delta}{3} \right\} \supset \bigcup\limits_{(t_j,t_j+\tau_0)} K(s),$$ 

\pa which satisfies

\begin{equation}\label{condOKBOUNDGEN}
d\left(\partial\Omega^+_{\delta}(t_j), \bigcup\limits_{(t_j,t_j+\tau_0)} K(s)\right) = \delta/3.
\end{equation}

\pa Notice that from (\ref{condAcotGen}) and (\ref{condOKBOUNDGEN}) we have

\begin{equation}\label{condOKBOUNDGEN2}
 d\left(\Omega^+_{\delta}(t_j), \bigcup\limits_{(t_j-\tau_0,t_j)} K(s)\right) > \delta/3.
\end{equation}

 \begin{figure}[H]
  \centering
    \includegraphics[width=8cm]{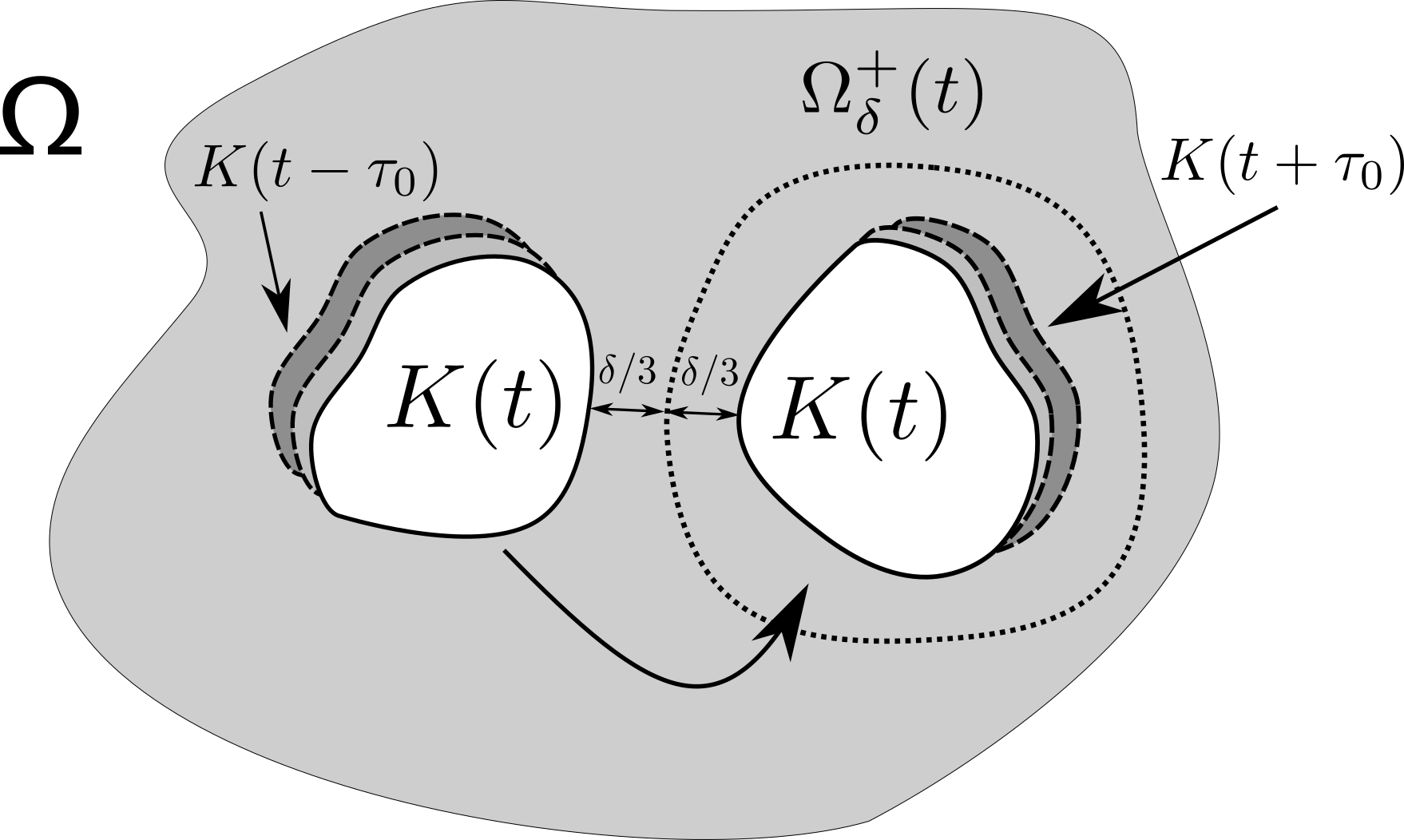}
        \caption{Neighbourhood $\Omega_\delta^+(t_j)$ stays at least at a distance $\delta/3$ from $\bigcup\limits_{s\in(t_j-\tau_0,t_j)} K(s)$ and its boundary stays at least at the same distance from $\bigcup\limits_{s\in(t_j,t_j+\tau_0)} K(s)$.}\label{figura}
\end{figure}

\pa Let $j\in \N$, since (\ref{condOKBOUNDGEN2}) holds, operating as in Theorem \ref{TeoAcot1} and applying Corollary \ref{acotlocSubdom} to the interval $(t_j-\tau_0,t_j)$ and the subset $\Omega^+_{\delta}(t_j)$ there exists a constant $\Lambda_1=\Lambda_1(\lambda,\rho,\tau_0,\nu,\delta)$ such that

$$u(t_j,\cdot) \leq \Lambda_1, \qquad \mbox{ for a.e. } x\in \Omega_{\delta}^+(t_j).$$

%Figura \ref{figura}

\pa Applying again Corollary \ref{acotlocSubdom} to the interval $(t_j-\tau_0,t_j+\tau_0)$ and the subset $\partial \Omega^+_{\delta}(t_j)$, %$\overline{\Omega\setminus\Omega^+_{\delta}(t_j)}$ 
there exists a constant $\Lambda_1^*=\Lambda_1^*(\lambda,\rho,\tau_0,\nu,\delta)$ such that

$$u(t,\cdot) \leq \Lambda_1^*, \qquad \forall t\in [t_j,t_j+\tau_0],\; x\in \partial\Omega_{\delta}^+(t_j).$$

\pa Observe that $\Lambda_1$ and $\Lambda_1^*$ do not depend on the initial condition nor the time $t_j$.  If we take $\Lambda = \max\{ \Lambda_1, \Lambda_1^*\}$,  which do not depend on the initial condition nor the time $t_j$ and we consider problem
\begin{equation}\left\{\begin{array}{lr}
v_t-\Delta v = \lambda v , & x\in\Omega_{\delta}^+(t_j), \; t\in (t_j,t_j+\tau_0),\\
v = \Lambda, & x\in\partial\Omega_{\delta}^+(t_j), \; t\in (t_j,t_j+\tau_0),\\
v(t_j,x)= \Lambda, & x\in\Omega_{\delta}^+(t_j),
\end{array}\right.\end{equation}
then, by comparison we have that $u(t,\cdot)\leq v(t,\cdot)$ for $t\in(t_j,t_j+\tau_0)$ and a.e. $x\in\Omega_{\delta}^+(t_j)$. 

Defining $w(t,x)=e^{-\lambda(t-t_j)}(v(t,x)-\Lambda)$, we have that $w$ satisfies 
\begin{equation}\left\{\begin{array}{lr}
w_t-\Delta w = \lambda \Lambda e^{-\lambda(t-t_j)} , & x\in\Omega_{\delta}^+(t_j), \; t\in (t_j,t_j+\tau_0),\\
w = 0, & x\in\partial\Omega_{\delta}^+(t_j), \; t\in (t_j,t_j+\tau_0),\\
w(t_j,x)= 0, & x\in\Omega_{\delta}^+(t_j).
\end{array}\right.\end{equation}
But, 
$$w(t,x)\leq \int_{t_j}^t\lambda\Lambda e^{-\lambda(s-t_j)}ds=\Lambda (1-e^{-\lambda(t-t_j)}),$$
which implies that
$$v(t,x)\leq \Lambda+e^{\lambda(t-t_j)}\Lambda (1-e^{-\lambda(t-t_j)})=\Lambda e^{\lambda(t-t_j)}.$$
Therefore, choosing $\Lambda^*=\Lambda e^{\lambda \tau_0}$ and taking into account that $u(t,x)\leq v(t,x)$, we prove the claim \eqref{claim}.

\medskip

\pa To show that the solution is bounded, we observe that by comparison
$$u(t,x)\leq U(t-(t_j+\tau_0),x), \qquad t\in [t_j+\tau_0,t_{j+1}+\tau_0],\quad j=1,2,\ldots$$
where $U(t,x)$ is the solution of the linear problem
\begin{equation}\label{eqUUU}
\left\{ \begin{array}{lr}
U_t - \Delta U = \lambda U , & x\in \Omega,\; t\in(0,\infty),\\
U=0, & x\in\partial\Omega, \; t\in(0,\infty),\\
U(0)=\Lambda^*, & x\in \Omega,
\end{array}\right.
\end{equation}
which obviously satisfies $U(t,x)\leq e^{\lambda t}\Lambda^*$ by the maximum principle and therefore 
$$\|u(t,\cdot)\|_{L^\infty(\Omega)}\leq \Lambda^*e^{|\lambda| \Xi}, \qquad t\geq t_1,$$
which shows the result.  \eproof

\rm 
Also, with a similar proof as Corollary \ref{corExtDomAcot1}, we can also show
\cor\label{corExtDomAcot4} Let $R>0$ such that $\overline{\Omega}\subset B(0,R)$. Then, assuming the same conditions as Theorem \ref{TeoAcot}, we have that all solutions of (\ref{PreEq1}) are bounded.

\rm

\medskip

\medskip
\section{Unboundedness}\label{secUnbound}
\rm
This section is devoted to obtaining conditions under which solutions grow up, that is, under which solutions become unbounded as $t\rightarrow \infty$.

\medskip

Observe that we have proven before, see Theorem \ref{TeoAcot1},  that if, roughly speaking, $\lambda<\lambda_0(K(t))$ $\forall t\geq t_0$, then we have boundedness of solutions. The equivalent situation for unboundedness would be to show that if $\lambda > \lambda_0(K(t))$ $\forall t\geq t_0$, the growth given by the linear part is strong enough to produce solutions which are unbounded as $t\rightarrow \infty$. But the previous section has also provided us with a clear counterexample to this hypothetical result. Actually we may have $K_0, K_1\neq \emptyset$ and $\lambda >\lambda_0(K_0)$, $\lambda>\lambda_0(K_1)$ but solutions remain bounded as $t\rightarrow \infty$ as detailed in Theorem \ref{TeoAcot} and seen in Example \ref{Kjump}. This again proves that the size of $\lambda_0(K(t))$ is not enough to decide the boundedness vs unboundedness of solutions. Actually,  the geometry and the speed at which the sets $K(t)$ move are also  important parameters.

\par\medskip 

 The assumptions in the following result can be interpreted as first of all having $K(t)$ large enough to allow the solution to grow inside these set. Moreover, the sets $K(t)$ move slowly enough so that the growth inside $K(t)$ is transported together with $K(t)$, allowing the function $u$ to grow without bounds inside the slowly moving set $K(t)$.

\par\medskip  We can state now:

\teo\label{TeoUnbounded} Let us assume that $n$ satisfies {\bf(N)}. Let us also assume that there exists a sequence of times $(t_i)_{i\geq 0}$ with $t_0<t_1<\ldots <t_i<\ldots$ and a family of nonempty, smooth  open subsets 

\begin{equation}\label{hypTeoUnbound} 
E_i\subset \bigcap\limits_{s\in [t_i,t_{i+1}]} K(s),
\end{equation}

\begin{figure}[H]
  \centering
    \includegraphics[width=8cm]{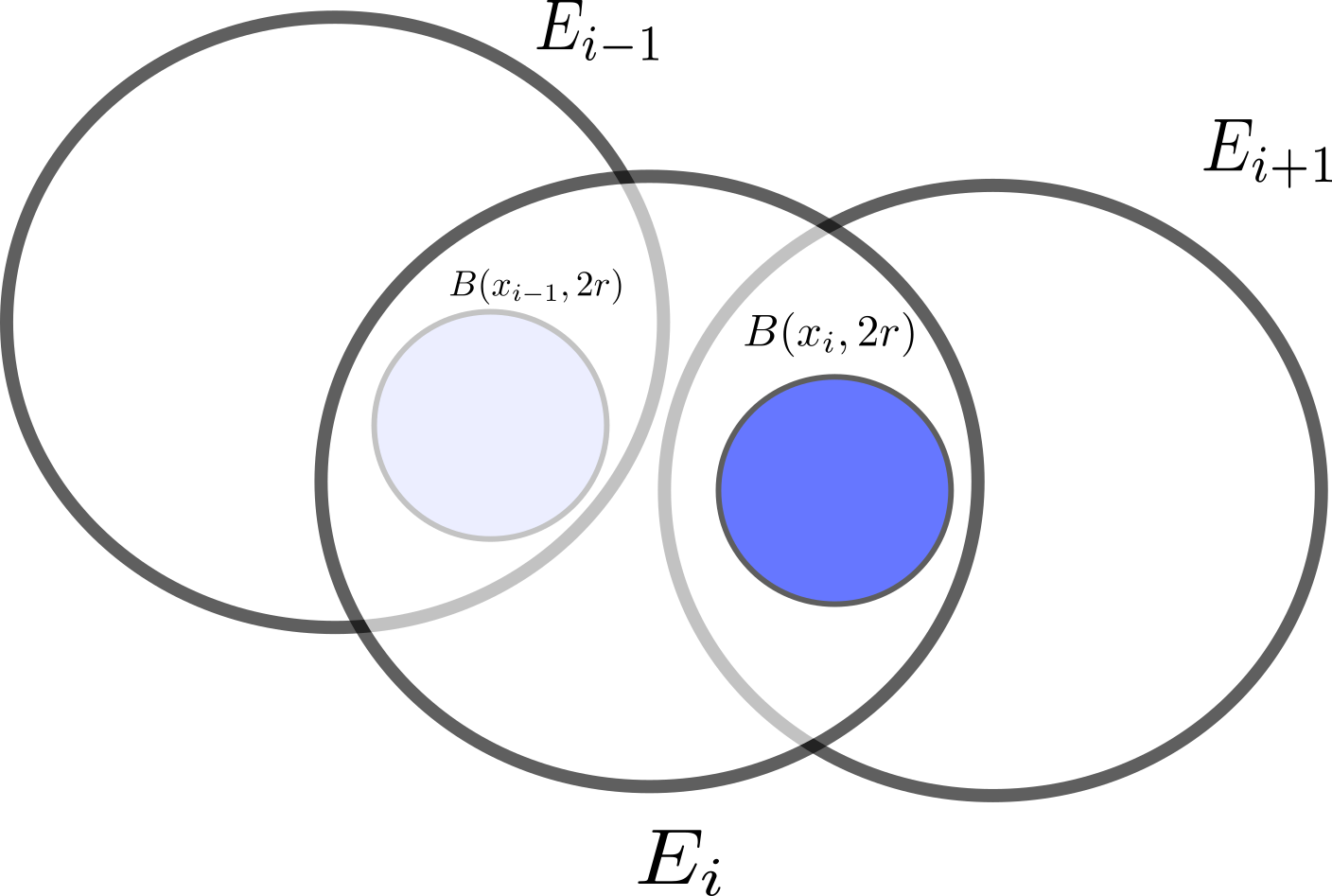}
        \caption{The sets $E_i$ and $B(x_i,2r)$}
        \label{figuraUnboundedEi2}
\end{figure}

\pa with rigid motions $T_i:\R^N\rightarrow \R^N$ such that $E_i := T_i(E_0)$, and such that there exists $(x_i)_{i\geq 0}$ and $r>0$  such that $B(x_i,2r)\subset E_i\cap E_{i+1}$ for $i=1,2,\ldots$, see Figure \ref{figuraUnboundedEi2}. Assume also that $\lambda>\lambda_1(E_0)$.

\pa Then, $\exists \tau>0$ such that if $t_{i+1}-t_i\geq \tau$ $\forall i=0,1,2,\ldots$ all solutions are unbounded.

\par\bigskip 
\rm Observe that having the sets $E_i$ satisfying  $\lambda>\lambda_0(E_i)$ is saying that the sets $E_i$ are large enough and therefore $K(t)$ are large enough from \eqref{hypTeoUnbound}. Also, intuitively, having  $t_{i+1}-t_i\geq \tau$ with property \eqref{hypTeoUnbound} means that the velocity at which the sets $K(t)$ move is slow enough.  

\par\medskip 

\rm 
The key point to prove this result is to find subsolutions which are unbounded as $t\rightarrow \infty$. Since the nonlinear logistic term is applied outside $K(t)$, unbounded growth will only be possible within $K(t)$. Therefore, we will consider the problem restricted to the sets $E_i$ which will be linear problems. Under certain conditions the solutions of these  linear problems will grow without bounds, providing us with the key factor to get the unboundedness for the nonlinear problem. 

\par\medskip  In this respect we have the following interesting auxiliary result on the growth of the linear problem and its relation with the first eigenfunctions.

\medskip

Let us denote by $\phi^A_n$ the n--th eigenfunction of the laplacian with homogeneous Dirichlet boundary conditions on the domain $A$ normalized in the $L^2$--norm and by $\lambda_n^A$ its eigenvalue. Assume also that $\phi_1^A$ is chosen to be positive in $A$. 

\lem\label{AuxUnbound} Let us consider the following linear problem \begin{equation}\left\{\begin{array}{ll}
v_t-\Delta v = \lambda v , & t>0,\; x\in E, \\
v = 0, & t> 0,\; x\in\partial E,\\
0\leq v(0,x)= v_0(x)\in L^\infty(E), & x\in E.
\end{array}\right.\end{equation}
where $E$ is a smooth bounded open connected domain. 

 \pa Let $D\subset E$ be a smooth open subset  with $\overline{D}\subset E$ and let  $\lambda>\lambda^E_{1}$. Then, if  $\gamma>1$ is given,  $v_0$ is not identically 0 and 
 $\tau>0$ is given by:
 \begin{equation}\label{def-tau}
 \tau=\max\left\{    \frac{N\lambda_2^E}{2e(\lambda_2^E-\lambda_1^E)}\left(\frac{2C_\infty\|v_0\|_{L^2(E)}}{\left\langle v_0, \phi^E_1 \right\rangle_{L^2(E)} \inf\limits_{D}{\phi^E_1}}\right)^{\frac{2}{N}},\frac{1}{\lambda-\lambda_1^E}\log \left(\frac{2\gamma\max\limits_{D}\phi_1^D}{\left\langle v_0, \phi^E_1 \right\rangle_{L^2(E)} \inf\limits_{D}{\phi^E_1} }\right)\right\},
 \end{equation}
 where $C_\infty=C_\infty(E,N)$ is the constant of the embedding $H^N(E)\hookrightarrow L^\infty(E)$, then 

$$v(t,x;v_0) \geq \gamma  \phi^{D}_1(x), \quad t\geq \tau,\; x\in D.$$

\proof  Let us define 
$$
I =\inf\limits_{D}{\phi^E_1}>0,\qquad
\mu=\frac{I}{\displaystyle 2\max\limits_{D} \phi_1^D},\qquad \alpha_j=\left\langle v_0, \phi^E_j \right\rangle_{L^2(E)},\quad j\in \N.$$

We have for $x\in D$ 
\begin{eqnarray*}
v(t,x)&=&\sum_{j=1}^{\infty}e^{\left(\lambda-\lambda^E_j\right)t} \langle v_0,\phi^E_j\rangle_{L^2(E)}\phi^E_j(x)\\
&=& e^{\left(\lambda-\lambda^{E}_1\right) t}\alpha_1 \phi^E_1(x)+\sum_{i=2}^{\infty}e^{\left(\lambda-\lambda^E_j\right)t} \alpha_j\phi^E_j(x)\\
&\geq& e^{\left(\lambda-\lambda^E_1\right)t}\alpha_1(\mu\phi^D_1+I/2)+\sum_{i=2}^{\infty}e^{\left(\lambda-\lambda^E_j\right)t} \alpha_j \phi^E_j(x)
\end{eqnarray*}
where we have used that
$$\phi_1^E(x)=\frac{\phi_1^E(x)}{2}+\frac{\phi_1^E(x)}{2}\geq \frac{I}{2}+\frac{I}{2}\geq \frac{I\phi_1^D}{2\max\limits_{D} \phi_1^D}+\frac{I}{2}=\mu\phi_1^D+\frac{I}{2}.$$

% Observe that $I>0$ since $E$ has smooth boundary,  $\phi^E_1>0$ and $\bar D\subset E$.

 Now, let us observe that

$$v(t,x)\geq e^{(\lambda-\lambda^E_1)t}\left(\alpha_1\mu\phi^D_1(x)+\alpha_1I/2 + J(x)\right),$$

\noindent where

$$J(x)=\sum_{j=2}^{\infty}e^{\left(\lambda^E_1-\lambda^E_j\right)t} \alpha_j\phi^E_j(x).$$

\noindent We want to arrive at an $L^\infty$--estimate for $J$. If $N$ is the dimension of the space and if $k>N/2$ we have the Sobolev inclusion $H^k(\Omega)\hookrightarrow L^\infty (\Omega)$. Thus, denoting by $C_{\infty}=C_{\infty}(N,E)$ the constant of the embedding $H^{N}(E)\hookrightarrow L^\infty(E)$,  we have 

\begin{eqnarray*}
\left\| J \right\|^2_{L^\infty(E)} & \leq & C_{\infty}^2\left\|\sum_{j=2}^{\infty}e^{\left(\lambda^E_1-\lambda^E_j\right)t} \alpha_j\phi^E_j\right\|^2_{H^{N}(E)} =  C_{\infty}^2\sum_{j=2}^{\infty}e^{-2\left(\lambda^E_j-\lambda^E_1\right)t} \alpha_j^2 (\lambda^E_j)^{N} \\
& = &C_{\infty}^2 \sum_{j=2}^{\infty}e^{-2\left(\lambda^E_j-\lambda^E_1\right)t} \left((\lambda^E_j-\lambda^E_1)t\right)^{N}\alpha_j^2 \frac{(\lambda^E_j)^{N}}{\left((\lambda^E_j-\lambda^E_1)t\right)^{N}} \\
& \leq  &C_{\infty}^2 \left(\frac{N}{2e}\right)^N \sum_{j=2}^{\infty} \alpha_j^2 \left(\frac{\lambda^E_j}{(\lambda^E_j-\lambda^E_1)t}\right)^N 
 \leq  C_{\infty}^2\left(\frac{N}{2e}\right)^N \sup_{j\geq 2}\left(\frac{\lambda^E_j}{(\lambda^E_j-\lambda^E_1)t}\right)^N \|v_0\|^2_{L^2(E)},
\end{eqnarray*}

\pa where we are bounding 
$$e^{-2\left(\lambda^E_j-\lambda^E_1\right)t} \left((\lambda^E_j-\lambda^E_1)t\right)^{N}\leq \sup_{a\geq 0}\{e^{-2a}a^N\}=\left(\frac{N}{2e}\right)^N, \quad \forall j\geq 2,$$
and using 
$$\sum_{j=2}^\infty\alpha_j^2\leq \sum_{j=1}^\infty\alpha_j^2=\|v_0\|^2_{L^2(E)}.$$

Moreover, since $\lambda_j^E\geq \lambda_2^E>\lambda_1^E$ for $j\geq 2$, we can easily check that  $\frac{\lambda^E_j}{\lambda^E_j-\lambda^E_1}$ is decreasig in $j\geq 2$ and therefore $\frac{\lambda^E_j}{\lambda^E_j-\lambda^E_1}\leq \frac{\lambda^E_2}{\lambda^E_2-\lambda^E_1}$. Hence,

\begin{eqnarray*} 
\left\| J \right\|^2_{L^\infty(E)}
& \leq &C_{\infty}^2\left(\frac{N}{2e}\right)^{N}\left(\frac{\lambda^E_2}{(\lambda_2^E-\lambda_1^E)t }\right)^{N}\|v_0\|_{L^2(E)}^2.
\end{eqnarray*}

\noindent Therefore, if we take   
$$t\geq    \frac{N\lambda_2^E}{2e(\lambda_2^E-\lambda_1^E)} \left(\frac{C_\infty 2\|v_0\|_{L^2(E)}}{\alpha_1 I}\right)^{2/N},$$
we have $\|J\|^2_{L^\infty(E)}\leq \alpha_1I/2$ and therefore, 

\begin{eqnarray*}
 v(t,x) \geq   e^{(\lambda-\lambda^E_1)t}\left(\alpha_1\mu\phi^{D}_1(x) +\alpha_1I/2 -\|J\|_{L^\infty(E)}\right)\geq 
 e^{(\lambda-\lambda^E_1)t}\alpha_1\mu\phi^{D}_1(x).
 \end{eqnarray*}

\pa Moreover, if $\gamma>1$ is given and we choose $t\geq \frac{1}{\lambda-\lambda_1^E}\log \left(\frac{\gamma}{\alpha_1\mu}\right)$, then  $e^{(\lambda-\lambda^E_1)t}\alpha_1\mu\geq \gamma$.  

Hence, if we take 

$$\tau=\max\left\{    \frac{N\lambda_2^E}{2e(\lambda_2^E-\lambda_1^E)} \left(\frac{C_\infty 2\|v_0\|_{L^2(E)}}{\alpha_1 I}\right)^{2/N}, \frac{1}{\lambda-\lambda_1^E}\log \left(\frac{\gamma}{\alpha_1\mu}\right)\right\}>0,$$

\pa we have  $v(t,x)\geq \gamma\phi_1^D(x)$ for all $t\geq \tau$ and $x\in D$.  This concludes the proof. \eproof

\par\bigskip

Let us consider now the following corollary, which will be important in the proof of the main result.

\cor\label{corolario-unbounded} If in the setting of the previous proposition, we assume that $B(a,2r)\subset E$ and take $v_0=\phi_1^{B(a,r)}$, then 
there exists $\rho>0$, independent of $a$, such that 

$$\|v(t,x,v_0)\|_{L^\infty(E)}\geq \rho, \quad t\geq 0.$$
\proof  Let us fix $a_0\in E$ such that $B(a_0,2r)\subset E$. Let us also fix $\gamma>1$.  Applying Proposition \ref{AuxUnbound}, with $D=B(a_0,r)\subset\subset E$, we have that the value of $\tau$ is independent of $a$. Hence, we have that 

$$v(t,x,\phi_1^{B(a,r)})\geq \gamma \phi_1^{B(a_0,r)}(x), \quad x\in B(a,r),\quad t\geq \tau.$$

This implies that for all $a\in E$  with $B(a,2r)\subset E$, we have 
$$\|v(t,\cdot, \phi_1^{B(a,r)})\|_{L^\infty(E)}\geq \gamma\| \phi_1^{B(a_0,r)}\|_{L^\infty(E)}=\gamma\| \phi_1^{B(0,r)}\|_{L^\infty(E)},  \quad \forall t\geq \tau. $$
Now, by comparison, with the problem posed in $B(a,r)$ and using that $\phi_1^{B(a,r)}$ is the first eigenfunction in $B(a,r)$, we have 
$$v(t,x,\phi_1^{B(a,r)})\geq e^{(\lambda-\lambda_1^{B(a,r)})t}\phi_1^{B(a,r)}(x), \quad t\geq 0,\quad x\in B(a,r).$$
But this implies that
$$\|v(t,\cdot,\phi_1^{B(a,r)})\|_{L^\infty(E)}\geq e^{-|\lambda-\lambda_1^{B(0,r)}|\tau}\|\phi_1^{B(0,r)}\|_{L^\infty(E)}, \qquad t\in [0,\tau].$$

Therefore, choosing
$$\rho=\min\{ \gamma\| \phi_1^{B(0,r)}\|_{L^\infty(E)}, e^{-|\lambda-\lambda_1^{B(0,r)}|\tau}\|\phi_1^{B(0,r)}\|_{L^\infty(E)}\},$$
we prove the result. \eproof

\par\medskip 

We are in a position now to give a proof of the main result, Theorem \ref{TeoUnbounded}.
\proof (of Theorem \ref{TeoUnbounded}). First let us state some immediate remarks derived from the assumptions. Since $B(x_i,2r)\subset E_i\cap E_{i+1}$ then $d(x_i,\partial E_i)\geq 2r$ and $d(x_i,\partial E_{i+1})\geq 2r$. Moreover, since $T_i:\R^N\rightarrow \R^N$ are rigid motions and $\lambda > \lambda_0(E_0)$ then $\lambda > \lambda_0(E_i)$ for any $i\geq 0$, and if we define $\tilde x_i=T_i^{-1}(x_i)$ then we have  $d(\tilde x_i,\partial E_0)\geq 2r$ for all $i=0,1,2,\ldots$.

\par\medskip To simplify, let us denote by $\psi$ the first eigenfuntion of $-\Delta$ in $B(0,r)$ with Dirichlet boundary conditions (which is explicitly known), that is $\psi=\phi_1^{B(0,r)}$. Also, denote by  $\psi_i:= \psi(x-x_i)$ the first eigenfunction of $-\Delta$ in $B(x_i,r)$ with Dirichlet homogeneous boundary conditions and extended by zero outside $B(x_i,r)$ for every $i\geq 0$.

\par\medskip We are going to apply Proposition \ref{AuxUnbound} for $i=1,2,\ldots$, with $E=E_i$, $D=B(x_i,r)$ and with initial condition $v_0=\psi_{i-1}$, which has support in $B(x_{i-1},r)\subset E_i$.  The value of $\tau$ obtained in the proposition, see (\ref{def-tau}),  is given for each $i$ by:

$$\tau_i=\max\left\{    \frac{N\lambda_2^{E_i}}{2e(\lambda_2^{E_i}-\lambda_1^{E_i})}\left(\frac{2C_{\infty}\|\psi_{i-1}\|_{L^2(E_i)}}{\left\langle \psi_{i-1}, \phi^{E_i}_1 \right\rangle_{L^2(E_i)} \inf\limits_{B(x_i,r)}{\phi^{E_i}_1}}\right)^{2/N}, \right.$$

$$\left.\frac{1}{\lambda-\lambda_1^{E_i}} \log \left(\frac{2\gamma\max\limits_{B(x_i,r)}\phi_1^{B(x_i,r)}}{\left\langle \psi_{i-1}, \phi^{E_i}_1 \right\rangle_{L^2(E_i)}\inf\limits_{B(x_i,r)}{\phi^{E_i}_1}} \right)\right\}>0.$$

Observe that since all sets $E_i$ are rigid transformations of the same set $E_0$,  then $\lambda_k^{E_i}=\lambda_k^{E_0}$, for all $k=1,2,..$, and the embedding constant $C_\infty$ is independent of $i$. Moreover since $B(x_i,r)\subset E_i^r\equiv \{x\in E_i: d(x,\partial E_i)>r\}$ and the maps $T_i$ are rigid motions, then   $\inf\limits_{B(x_i,r)}\phi_1^{E_i}\geq \inf\limits_{E_i^r}\phi_1^{E_i}=\inf\limits_{E_0^r}\phi_1^{E_0}$, which is independent of $i$.  Moreover, $\|\psi_{i-1}\|_{L^2(E_i)}=1$ and $\max\limits_{B(x_i,r)}\phi_1^{B(x_i,r)}=\max\limits_{B(0,r)}\phi_1^{B(0,r)}>0$ which are also independent of $i$. Also,

$$\left\langle \psi_{i-1}, \phi^{E_i}_1 \right\rangle_{L^2(E_i)}=\int_{E_i}\psi_{i-1}\phi_1^{E_i}\geq \left(\int_{B(x_{i-1},r)}\psi_{i-1}\right)\inf\limits_{B(x_{i-1},r)}\phi_{1}^{E_i}$$
$$\geq \left(\int_{B(0,r)}\psi\right)\inf\limits_{E_i^r}\phi_{1}^{E_i}=\left(\int_{B(0,r)}\psi\right)\inf\limits_{E_i^0}\phi_{1}^{E_0}>0,$$

\pa which is also independent of $i$.  Therefore, if we define 

$$\tau=\max\left\{   \frac{N\lambda_2^{E_0}}{2e(\lambda_2^{E_0}-\lambda_1^{E_0})}\left(\frac{2C_\infty}{\int_{B(0,r)}\psi (\inf\limits_{E_0^r}{\phi^{E_0}_1})^2}\right)^{2/N}, \frac{1}{\lambda-\lambda_1^{E_0}} \log \left(\frac{2\gamma\max\limits_{B(0,r)}\psi}{\int_{B(0,r)}\psi (\inf\limits_{E_0^r}{\phi^{E_0}_1})^2} \right)\right\}>0,$$

\pa we can apply the previous proposition for intervals of time larger or equal than $\tau$.

\par\bigskip

Hence, if for $i=1,2\ldots$, we have  $t_{i}-t_{i-1}\geq \tau$ and we consider  problems
\begin{equation}\label{EqUnCircling}\left\{\begin{array}{ll}
v_{i_t}-\Delta v_i = \lambda v_i , & x\in E_{i}, \; t\in (t_{i-1},t_{i}],\\
v_i(t,x) = 0, & x\in\partial E_{i}, \; t\in (t_{i-1},t_{i}],\\
v_i(t_{i-1},x)=\gamma^{i-1}\psi_{i-1}, & x\in E_{i},
\end{array}\right.\end{equation}

\pa then we have that $v_i(t_i,x)\geq \gamma^i\psi_{i}(x)$ in $B(x_{i+1},r)$.  If we extend the functions $v_i(t,\cdot)$ by 0 to all $\Omega$ and define the function $v(t,x)$ as  $v(t,x)=v_i(t,x)$ $t\in [t_{i-1},t_i)$, then by comparison it is not difficult to see that 
$$v(t,x)\leq u(t,x,t_0,\psi_0(\cdot)),\quad t\geq t_0.$$

Indeed, if $t\in [t_0,t_1]$ we have that $v_1(t,x)$ is a subsolution of $u(t,x)$ and therefore $u(t_1,x)\geq v_1(t_1,x)\geq \gamma \psi_1(x)$.  Observing that $\gamma \psi_1(x)$ is the initial condition of $v_2$, we get again that in $[t_1,t_2]$ the function $v_2(t,x)$ is a subsolution of $u(t,x)$.  Therefore $u(t_2,x)\geq v_2(t_2,x)\geq \gamma^2 \psi_2(x)$. By induction it is not difficult to see that in $[t_{i-1},t_i]$ we have that $v_i(t,x)$ is a subsolution of $u(t,x)$ and therefore $u(t_i,x)\geq \gamma^i\psi_i(x)$.  This means that 

$$u(t,x)\geq v(t,x),\qquad t\geq t_0, \quad x\in \Omega,$$
and in particular
$$\|u(t_i,x)\|_{L^\infty(\Omega)}\geq \gamma^i\|\psi_i\|_{L^\infty(\Omega)}=\gamma^i \|\psi\|_{L^\infty(B(0,r))}\to \infty, \hbox{ as } i\to \infty.$$

Moreover, with Corollary \ref{corolario-unbounded} we can actually see that  $\|v_i(t,x)\|_{L^\infty(\Omega)}\geq \gamma^{i-1}\rho$ for each $t\in [t_{i-1},t_i]$. Hence, we get
$$\lim_{t\to+\infty} \|u(t,\cdot)\|_{L^\infty(\Omega)}=+\infty.$$ 
This concludes the proof of the theorem. \eproof

\par\medskip

\re Observe that if we have the situation where $E_i\equiv E_0$ for all $i=1,2,\ldots$, condition \eqref{hypTeoUnbound}  is read as $E_0\subset K(t)$ for all $t\geq t_0$ and since we are assuming  $\lambda>\lambda_0(E_0)$ we can show that the solutions are unbounded. For this, notice that $n(t,x)\leq \tilde n(x)$ where $\tilde n(x)=\|n\|_{L^\infty}\chi_{\Omega\setminus E_0}(x)$. Via comparison with the autonomous problems with $\tilde n$, we show that the solutions are all unbounded. \rm

\par\medskip 

In the line of the above remark, in the next result we show that the assumption that $E_0\subset K(t)$ for all $t$ can be weakened.  For this, let us consider problem (\ref{PreEq1}) where $n$ satisfy {\bf (N)} and $K(t)$ satisfies

\begin{equation}\label{subn0}\left\{\begin{array}{lr}
K(t)\supset K_1, & t\in(T_{0,i},T_{1,i}], x\in \Omega,\\
K(t)\supset K_0, & t\in(T_{1,i},T_{0,i+1}], x\in \Omega,
\end{array}\right.\end{equation}

\pa for some unbounded time sequences $t_0\leq T_{0,i}\leq T_{1,i}\leq T_{0,i+1} \rightarrow \infty$ as $i\rightarrow \infty$ and for some sets $K_0=\bar{\Omega_0}$, $K_1=\bar{\Omega_1}$ with $\Omega_0\subset \Omega_1$, smooth open sets.

\par\medskip We can show the following 

\prop\label{PropUnbound}
With the notation above, if we assume 
\begin{equation}\label{lambdapos}
\lambda_0(K_1)<\lambda<\lambda_0(K_0)<\infty.
\end{equation}

\pa and that there exist $\tau$ and $\eta$  such that for the time sequences  $t_0\leq  T_{0,i} < T_{1,i} <T_{0,i+1}\rightarrow \infty$ satisfy 

\begin{equation}\label{seqTB}\left\{\begin{array}{lr}
\left|T_{1,i}-T_{0,i}\right|> \tau, & i\in\N,\\
\left|T_{0,i+1}-T_{1,i}\right| \leq \eta, & i\in\N,
\end{array}\right.\end{equation}

\pa then, if $\tau$ is sufficiently large, the solution $u(t,x)$  of (\ref{PreEq1}) is unbounded. 

% $$\limsup_{t\to \infty}\|u(t;u_0)\|_{L^{\infty}(\Omega)}=\infty.$$

\proof
We want to exhibit a subsolution of $u(t,x)$ that grows up as time goes to infinity. Without any loss of generaliy, we can assume that 

\begin{equation}\label{condlanda}
\lambda \in (\lambda^{\Omega_1}_1, \lambda^{\Omega_1}_2),
\end{equation}

\pa since, if greater or equal than $\lambda^{\Omega_1}_2$, we can take a smaller $\lambda$ which verifies (\ref{condlanda}) and it would result in a subsolution of the original solution.  Let us consider the following family of problems for each $i\geq 1$:

\begin{equation}\label{sub1}
\left\{\begin{array}{lll}
U_{i_t}(t,x)-\Delta U_i(t,x) = \lambda U_i(t,x),\ \ &t\in(T_{0,i},T_{1,i}],& x\in \Omega_1,\\
U_i(t,x) = 0,&t\in(T_{0,1},T_{1,i}],& x\in \partial \Omega_1,\\
U_i(T_{0,i},x) = \phi_1^{\Omega_0}, &  & x\in \Omega_1,
\end{array}\right.
\end{equation}

\begin{equation}\label{sub2}
\left\{\begin{array}{lll}
W_{i_t}(t,x)-\Delta W_i(t,x) = \lambda W_i(t,x),\ \ &t\in(T_{1,i},T_{0,i+1}],& x\in \Omega_0,\\
W_i(t,x) = 0,&t\in(T_{1,i},T_{0,i+1}],& x\in \partial \Omega_0,\\
W_i(T_{1,i},x)=\phi_1^{\Omega_0}, & & x\in \Omega_0.
\end{array}\right.
\end{equation}

Observe that $U_i$ is posed in the large set $\Omega_1$ and therefore the solution there tends to grow, while $W_i$ is posed in the smaller set $\Omega_0$ where the solution tends to decay.  As a matter of fact,  
$$W_i(t)=e^{(\lambda-\lambda_1^{\Omega_0})(t-T_{1,i})}\phi_1^{\Omega_0},$$
and if we define  
$$\alpha= e^{(\lambda-\lambda_1^{\Omega_0})\eta}\in (0,1),$$
we have
\begin{equation}\label{estimate2}
  W_i(t)\geq \alpha \phi_1^{\Omega_0}, \qquad t\in [T_{1,i},T_{0,i+1}].
  \end{equation}
%where the constants $b_i$ will be chosen appropriately. 

\par\medskip  Applying now Proposition \ref{AuxUnbound}  to $U_i$, with $E=\Omega_1$, $D=\Omega_0$, and $U_{i,0}=\phi_1^{\Omega_0}$, we have that if $\gamma>1$ and define 

 $$ \tau=\max\{    \frac{N\lambda_2^{\Omega_1}}{2e(\lambda_2^{\Omega_1}-\lambda_1^{\Omega_1})}\left(\frac{2C_\infty\|\phi_1^{\Omega_0}\|_{L^2(\Omega_1)}}{\left\langle \phi_1^{\Omega_0}, \phi^{\Omega_1}_1 \right\rangle_{L^2(\Omega_1)} 
 \inf\limits_{\Omega_0}{\phi^{\Omega_1}_1}}\right)^{2/N},\frac{1}{\lambda-\lambda_1^{\Omega_1}}
 \log \left(\frac{2\frac{\gamma}{\alpha}\max\limits_{\Omega_0}\phi_1^{\Omega_0}}{\left\langle \phi_1^{\Omega_0}, \phi^{\Omega_1}_1 \right\rangle_{L^2(\Omega_1)} \inf\limits_{\Omega_0}{\phi^{\Omega_1}_1} }\right)\},$$
which does not depend on $i$.  Then,  choosing  $T_{1,i}-T_{0,i}\geq \tau$, we have 
\begin{equation}\label{estimate-Ui-1}
U_i(t)\geq \frac{\gamma}{\alpha} \phi_1^{\Omega_0}, \quad t\in [T_{0,i}+\tau, T_{1,i}].
\end{equation}

Also, notice that comparing problem \eqref{sub1}) with the same problem but posed in $\Omega_0\subset \Omega_1$, we have that
\begin{equation}\label{estimate-Ui-2}
U_i(t)\geq e^{(\lambda-\lambda_1^{\Omega_0})t}\phi_1^{\Omega_0}\geq e^{(\lambda-\lambda_1^{\Omega_0})\tau}\phi_1^{\Omega_0}, \quad t\in [T_{0,i}, T_{0,i}+\tau].
\end{equation}
Putting together \eqref{estimate-Ui-1} and \eqref{estimate-Ui-2}, we get 
\begin{equation}\label{estimate1}
U_i(t)\geq \rho\phi_1^{\Omega_0},\qquad t\in [T_{0,i},T_{1,i}]
\end{equation}
where $\rho=\min\{\frac{\gamma}{\alpha},e^{(\lambda-\lambda_1^{\Omega_0})\tau}\}$.
\par\medskip Let us construct the function $Z(t,x)$ where for every $i\in\N$:

\begin{equation*}Z(t,\cdot)=\left\{\begin{array}{lr}
 \gamma^{i-1} U_i(t,\cdot), & t\in (T_{0,i},T_{1,i}],\\
\frac{\gamma^{i}}{\alpha}W_i(t,\cdot), & t\in (T_{1,i},T_{0,i+1}],\\
\end{array}\right.\end{equation*}

\pa where we assume $ U_i, W_i$ are extended by 0 to the whole domain $\Omega$.  With the analysis above, we can conclude that $Z(t,x)$ is a subsolution of $u(t,x)$ and with \eqref{estimate2} and \eqref{estimate1} we get that $\lim_{t\to+\infty}\|Z(t,\cdot)\|_{L^\infty}=+\infty$, which proves the result. \eproof

\par\bigskip\noindent {\bf \large Declarations.}

\par\medskip\noindent{ \bf Funding.}   {\bf José M. Arrieta} is partially supported by grants PID2019-103860GB-I00, PID2022-137074NB-I00  and CEX2019-000904-S ``Severo Ochoa Programme for Centres of Excellence in R\&D'' the three of them from MICINN, Spain. Also by ``Grupo de Investigación 920894 - CADEDIF'', UCM, Spain.   {\bf Marcos Molina-Rodríguez} is partially supported by the predoctoral appointment  BES-2013-066013 and grant PID2019-103860GB-I00 both by MICINN, Spain and by ``Grupo de Investigación 920894 - CADEDIF'', UCM, Spain.  {\bf Lucas A. Santos} is partially supported by the Coordenação de Aperfeiçoamento de Pessoal de Nível Superior – Brasil (CAPES) – Finance Code 001”

\par\medskip\noindent {\bf Ethical Approval }.  Not applicable

\par\medskip\noindent{\bf Availability of data and materials} Not applicable.

\end{document}